\documentclass[10pt]{article}

\usepackage{authblk}

\usepackage{amsmath, amssymb}
\usepackage{stmaryrd}               % for llbracket/rrbracket
\usepackage{bm}                     % for bold math font

\usepackage{libertine}		
\usepackage{libertinust1math}
\usepackage{graphicx}
\usepackage{fancyhdr}
\usepackage{tcolorbox}              % frame

\usepackage[ddmmyyyy]{datetime}           
\usepackage{microtype}              % fixes problems with "Overfull \hbox"
\usepackage{float}                  % fix position of figures
\usepackage[hidelinks]{hyperref}    % clickable links and references

\newcommand{\R}{\mathbb{R}}
\newcommand{\F}{\mathcal{F}}
\newcommand{\T}{\mathcal{T}}
\newcommand{\M}{\mathcal{M}}
\newcommand{\dt}{\text{dt}}
\newcommand{\patch}{\mathcal{P}_F}
\newcommand{\weight}{\omega}

\title{A High-Order Conservative Cut Finite Element Method for Problems in Time-Dependent Domains}
\author[a]{Sebastian Myrbäck\thanks{Corresponding author:\\ \textit{Email addresses:} smyrback@kth.se (Sebastian Myrbäck), sara.zahedi@math.kth.se (Sara Zahedi)}}
\author[a]{Sara Zahedi}

\affil[a]{Department of Mathematics, KTH Royal Institute of Technology, SE-100\,44 Stockholm, Sweden.}

\date{}

\begin{document}

\maketitle

\begin{abstract}
    \noindent A mass-conservative high-order unfitted finite element method for convection-diffusion equations in evolving domains is proposed. The space-time method presented in [P. Hansbo, M. G. Larson, S. Zahedi, Comput. Methods Appl. Mech. Engrg. 307 (2016)] is extended to naturally achieve mass conservation by utilizing Reynold's transport theorem. Furthermore, by partitioning the time-dependent domain into macroelements, a more efficient stabilization procedure for the cut finite element method in time-dependent domains is presented. Numerical experiments illustrate that the method fulfills mass conservation, attains high-order convergence, and the condition number of the resulting system matrix is controlled while sparsity is increased. Problems in bulk domains as well as coupled bulk-surface problems are considered.
\end{abstract}

\noindent {\bf Keywords:} Unfitted finite element method, space-time finite element method, convection-diffusion equation, Reynold's transport theorem, mass conservation, soluble surfactant

\maketitle

%% Paper
\section{Introduction}\label{section:introduction}

Unfitted, immersed, or embedded finite element methods have become a popular alternative to standard finite element methods for solving partial differential equations (PDEs) in domains with complex or evolving boundaries~\cite{cutfem-overview, unfitted-methods, finite-cell, ext_fem, fictfem1, fictfem2, Olshanskii}.

We are interested in high-order unfitted discretizations of convection-diffusion equations in time-dependent domains. The motivation is the study of so-called surface-active agents, or surfactants, which are substances that reduce the surface tension of a liquid. They play a crucial role in numerous natural and industrial processes, such as the stabilization of emulsions in food and cosmetic products, enhanced oil recovery, and biological processes such as in lung function~\cite{surfactants, lungs, surf_appl}.

In~\cite{hansbo2016cut, Zah17} a space-time cut finite element method (CutFEM) is proposed for solving convection-diffusion equations in time-dependent domains. The method avoids remeshing procedures as the domain is moving.  Other unfitted discretizations including high-order methods for convection-diffusion equations in bulk domains have been proposed, see~\cite{lehrenfeld, trackers, BADIA202360, HaLaZa15}. However, none of these discretizations conserve mass. In~\cite{hansbo2016cut, Zah17}, mass conservation is achieved by enforcing it through a Lagrange multiplier. Therefore, in~\cite{frachon2022cut, SeMy22} we started studying how mass conservation can be inherited more naturally in unfitted discretizations by reformulating the weak form using Reynolds' transport theorem. 

The strength of the cut finite element method in~\cite{hansbo2016cut, Zah17} is its simple implementation, which is made possible by the ghost-penalty stabilization added in the method. In this method, the space-time domain is never explicitly constructed and space-time integrals are split into a weighted sum of space-integrals evaluated at discrete quadrature points in time.

In a cut finite element method, the domain of interest is embedded in a computational domain equipped with a background mesh that is unfitted with respect to the boundary of the domain. An active mesh is defined, on which appropriate finite element spaces and a weak form are formulated. In the space-time CutFEM, the active mesh consists of elements in the background mesh that have an intersection with the time-dependent domain in a given time interval. The finite element spaces utilize discontinuous elements in time. The weak formulation of the PDE includes ghost penalty stabilization terms. The ghost penalty stabilization technique first introduced in~\cite{ghost-penalty-stabilization} yields a well-defined extension of the approximate solution to the entire active mesh. In general, a strategy for handling elements that are badly cut is necessary for all unfitted methods since these elements otherwise cause problems with singular or near-singular linear systems. Adding ghost-penalty stabilization has become the most common strategy in connection with cut finite element methods see e.g.~\cite{cutfem-overview,hansbo2016cut, surface-cutfem-2019, burman2015, LehOls19}. 

Recently, a strategy on how to apply ghost-penalty only where necessary was developed for stationary domains in~\cite{larson2021conservative}. The idea is to construct a macroelement partition of the mesh, where macroelements have a large intersection with the domain of interest, and therefore stabilization is only needed on the internal faces. In this paper, we develop a time-dependent macroelement stabilization, drawing upon ideas from~\cite{larson2021conservative}. We are able to lower the $L^2$ error and reduce its sensitivity to stabilization parameters, ensure control of the condition number of the system matrix, and increase the sparsity of the system matrix compared to when full stabilization is used.

The discretization in this paper achieves mass conservation naturally for both bulk and coupled bulk-surface convection-diffusion equations by using Reynold's transport theorem. Furthermore, to accurately compute the integrals on elements cut by the interface we use quadrature rules from~\cite{saye}. For a different strategy, using isoparametric mappings, see~\cite{lehrenfeld, lehrenfeld_stationary}. However, obtaining mass conservation seems to be challenging with the isoparametric method, see Remark 19 in~\cite{preuss2018}. The stability of the conservative method presented here is ensured by utilizing macroelement stabilization for time-dependent domains and using large stabilization where it is needed, without introducing larger errors. We are not aware of any other unfitted discretization that achieves both mass conservation and is higher order than second order.  

The paper is outlined as follows. In Section~\ref{section:model-problem}, we present the mathematical model, a convection-diffusion equation in a time-dependent domain. In Section~\ref{section:mesh_fespace}, we define the mesh and the finite element spaces of the space-time CutFEM.\@ We present the conservative discretization in Section~\ref{section:convection-diffusion}. In Section~\ref{section:stabilization}, we present a novel macroelement stabilization in time-dependent domains. In Section~\ref{section:numerical_quadrature}, we present the numerical quadrature rules used in the method. Various numerical examples illustrating the performance of the method can be found in Section~\ref{section:numerical-examples}. In Section~\ref{section:coupled_bulk_surface}, we extend the ideas to the coupled bulk-surface problem and show numerical results. Finally, we conclude our results in Section~\ref{section:conclusion}.

\section{The Mathematical Model}\label{section:model-problem}

\begin{figure}
    \centering
    \includegraphics[width=.45\textwidth]{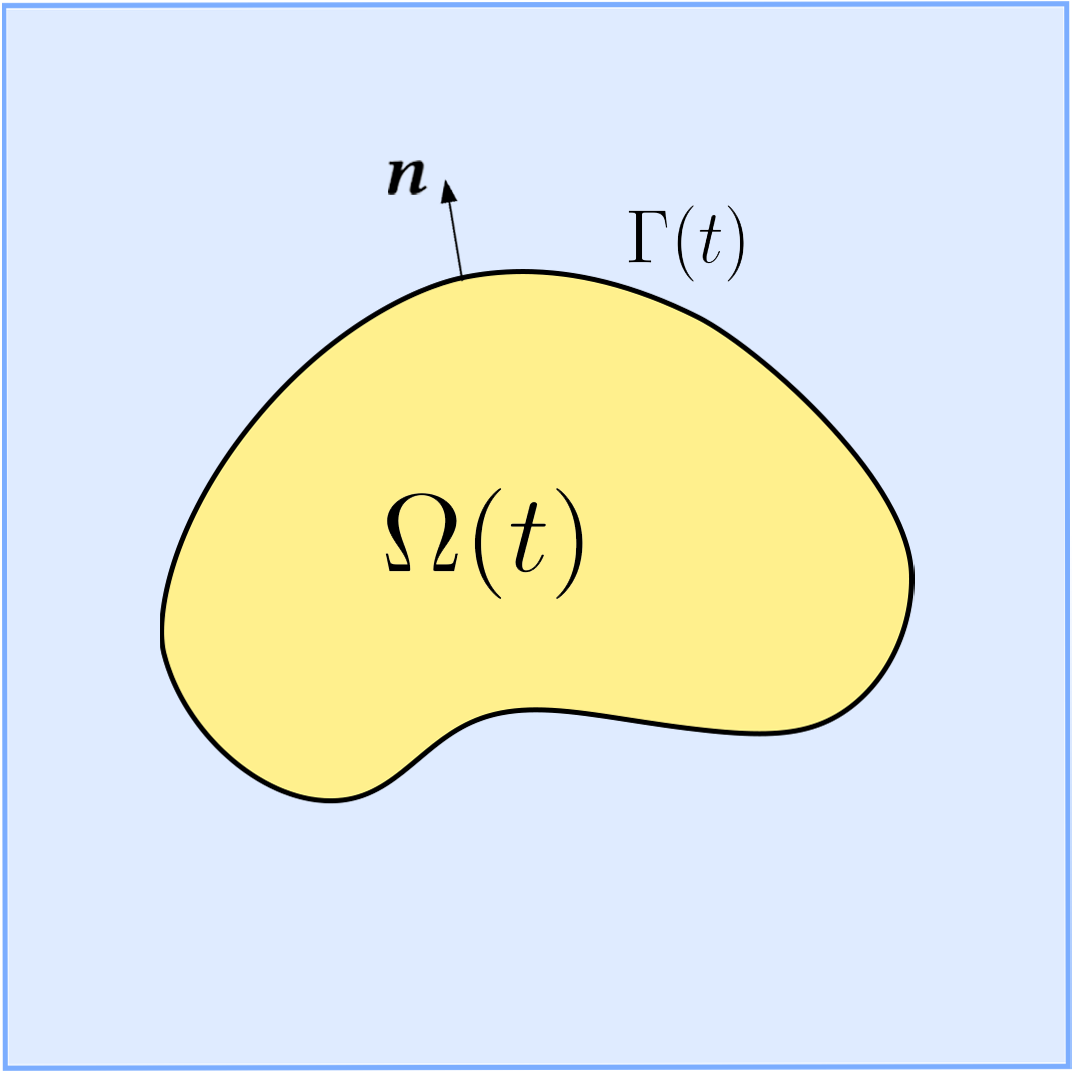}
    \caption{A domain sketch.}\label{fig:domain}
\end{figure}

For $t\in I=[0,T]$ with $0 < T < \infty$, let $\Omega(t) \subset \R^d$, $d=2,3$, denote a bounded domain with boundary $\partial\Omega(t) = \Gamma(t) \subset \R^{d-1}$, and let $\Omega_0$ be a polytopal open subset of $\R^d$ such that for every $t\in I$, $\Omega(t)\subset\Omega_0$. See Figure~\ref{fig:domain} for an illustration.
The surface $\Gamma(t)$ evolves with a given, sufficiently smooth, velocity field $\bm{\beta}: I \times \Omega_0 \rightarrow \R^d$ that is assumed to be divergence-free. Let $\bm{n}$ denote the unit normal vector on $\Gamma(t)$, pointing outward from $\Omega(t)$. 

\subsection{The convection-diffusion equation}\label{subsection:convection-diffusion-bulk}

We consider the following convection-diffusion equation in the evolving domain $\Omega(t)$: Find $u: I\times \Omega(t) \rightarrow \R$ such that 
\begin{alignat}{2}
    \partial_t u + \nabla\cdot(\bm{\beta} u) - \nabla\cdot\left(D\nabla u\right) &= f, \qquad\quad&&\text{in } \Omega(t),\label{eq:timebulk}\\
    \bm{n}\cdot D \nabla u &= 0, &&\text{on }\Gamma(t),\label{eq:neumann}\\
    u(0, \bm{x}) &= u_0(\bm{x}), &&\text{in }\Omega(0)\label{eq:initial-condition-bulk},
\end{alignat}
where $t\in I$, $f: I\times \Omega(t)\rightarrow \R$ is a source/sink function, and $u_0: \Omega(0)\rightarrow \R$ is the known initial solution. 

Recall Reynold's transport theorem
\begin{align}
    \frac{\text{d}}{\dt}\int_{\Omega(t)}u &= \int_{\Omega(t)} \partial_t u + \int_{\Gamma(t)}\bm{n}\cdot\bm{\beta} u.\label{eq:reynold-bulk}
\end{align}
Integrating~\eqref{eq:timebulk} over $\Omega(t)$, and using the divergence theorem yields
\begin{equation}
    \int_{\Omega(t)}\partial_t u + \int_{\Gamma(t)}\bm{n}\cdot\bm{\beta}u - \int_{\Gamma(t)}\bm{n}\cdot D\nabla u = \int_{\Omega(t)}f .\label{eq:integrated}
\end{equation}
Combining~\eqref{eq:reynold-bulk} and~\eqref{eq:integrated}, and using the Neumann boundary condition~\eqref{eq:neumann} yields
\begin{align}
    \frac{\text{d}}{\dt} \int_{\Omega(t)} u = \int_{\Omega(t)}f\label{eq:conservation-law-bulk}.
\end{align}
Integrating this equation over $I = [0,T]$ we get
\begin{equation}
    \int_{\Omega(T)}u - \int_{\Omega(0)} u = \int_0^T\int_{\Omega(t)} f,
    \label{eq:conservation-law-non-zero}
\end{equation}
which asserts that the change of mass is equal to the amount of mass produced inside the domain. In Section~\ref{section:convection-diffusion}, we present a numerical scheme where the solution fulfills this relation.

\section{Mesh and Finite Element Spaces}\label{section:mesh_fespace}

In this section, we define the active mesh and the finite element spaces for the space-time CutFEM~\cite{hansbo2016cut}. 

\subsection{Mesh}
Let $\{\T_h\}_h$ be a quasi-uniform family of meshes of the domain $\Omega_0$, with mesh parameter $0<h\leq h_0<<1$ such that $\Omega_0 = \bigcup_{K\in\T_h}K$. We will use quadrilateral elements. Let $0=t_0 < t_1 < \cdots < t_N = T$ be a partition of $I = [0,T]$ into time intervals $I_n=(t_{n-1},t_n]$ of length $\Delta t_n = t_{n}-t_{n-1}$ for $n=1,\dots,N$.

\begin{figure}
    \centering
    \includegraphics[width=0.475\textwidth]{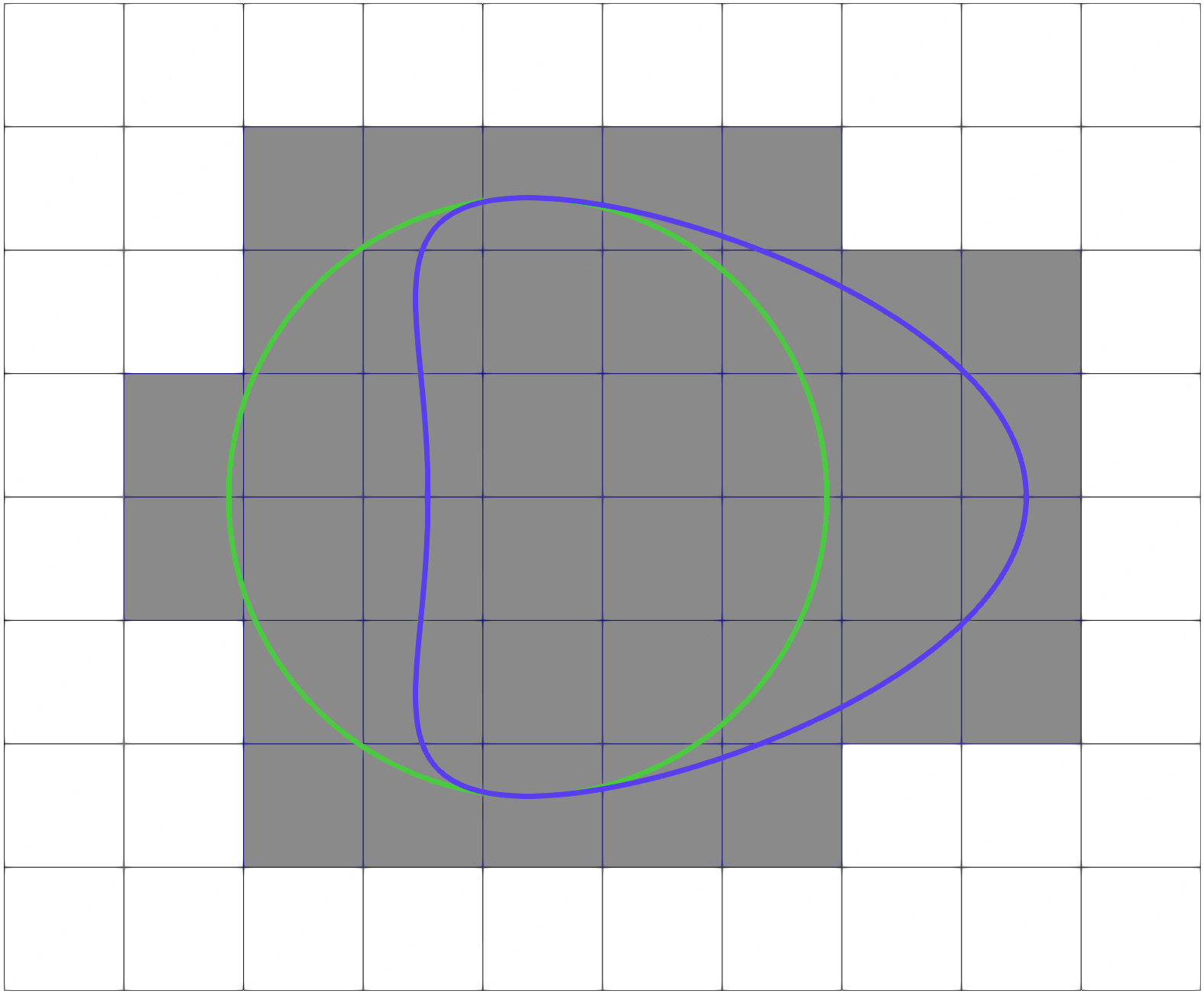}%
    \hfill
    \includegraphics[width=0.475\textwidth]{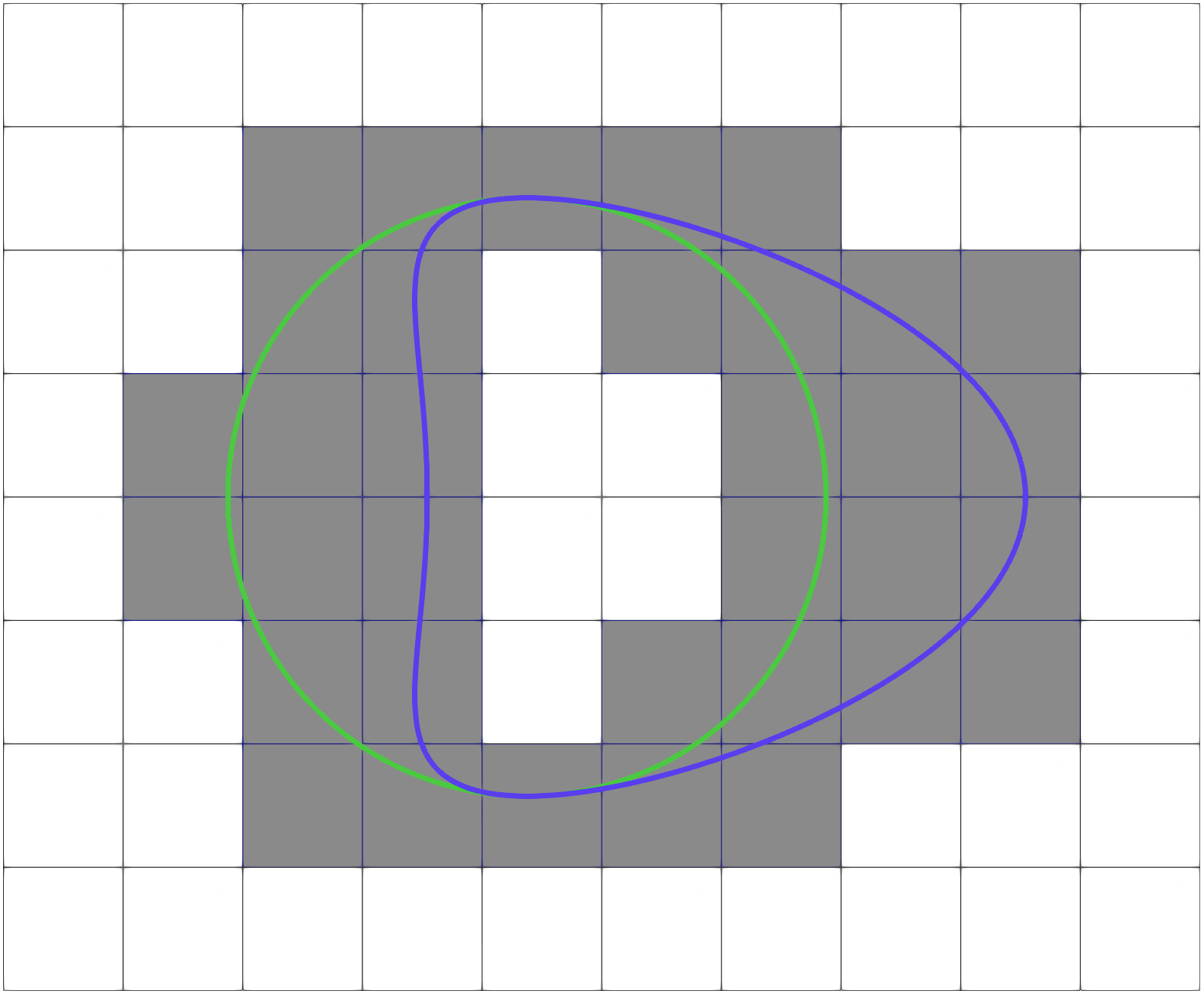}
    \caption{Background domain in white and elements in the active meshes $\T_{h}^n$ (left) and $\T_{h,\Gamma}^n$ (right) in grey. The boundary $\Gamma(t)$ is visualized at time instances $t=t_{n-1}$ (green) and $t=t_n$ (blue).}\label{fig:active_meshes}
\end{figure}

For each time interval $I_n$ we define active meshes associated to $\Omega(t)$ and to the boundary $\Gamma(t)$ as:
\begin{align}
    \T_h^n &= \{K \in \T_h : K\cap \Omega(t) \neq \emptyset \text{ for some } t\in I_n\},\label{eq:active-mesh-time-bulk}\\
    \T_{h,\Gamma}^n &= \{K \in \T_h : K\cap \Gamma(t) \neq \emptyset \text{ for some } t\in I_n\},\label{eq:active-mesh-time-surf}
\end{align}
These active meshes constitute the following active domains
\begin{align}
    \Omega_{\T_h^n} &= \bigcup_{K\in \T_{h}^n} \{K\},\\ 
    \Omega_{\T_{h,\Gamma}^n} &= \bigcup_{K\in \T_{h,\Gamma}^n} \{K\}.
\end{align}
For a visualization of the active meshes, see Figure~\ref{fig:active_meshes}. 

\subsection{Finite element spaces}

Let $P_k(I_n)$ denote the space of polynomials of degree $k$ in $I_n$, and let $V_{h,m}$ be the tensor product space of Lagrange basis functions of order $m$ in each space direction on the quadrilateral elements on the background mesh $\T_h$. We define the space-time cut finite element space as 
\begin{equation}
    W_{h,k,m}^n = P_k(I_n)\otimes V_{h,m}|_{\Omega_{\T_h^n}}.
\end{equation}
Due to the tensor product structure of $W_{h,k,m}^n$, functions in this space are on the form
\begin{align}
    v(t,\bm{x}) &= \sum_{i=0}^k v_{h, k}(\bm{x})\left(\frac{t-t_{n-1}}{\Delta t_n}\right)^k,\quad t\in I_n,\ \bm{x}\in \Omega_{\T_{h}^n},
\label{eq:time-functions}
\end{align}
where
\begin{align}
    v_{h,k}(\bm{x}) = \sum_{l=1}^{N_{m}^n}c_{l,k}\varphi_l(\bm{x}),
\end{align}
where $c_{k,l}\in \R$ and $N_{m}^n$ is the number of degrees of freedom of $V_{h,m}|_{\Omega_{\T_h^n}}$ for which $\{\varphi_l\}_l$ is a basis. 

\section{Space-Time Cut Finite Element Discretizations}\label{section:convection-diffusion}
In this section, we state two space-time cut finite element discretizations for the convection-diffusion equation.  

\subsection{Conservative scheme}\label{subsection:conservative-scheme}
One can apply Reynold's transport theorem~\eqref{eq:reynold-bulk} to the product $uv$, where $u,v\in H^1(\Omega(t))$ (see e.g. Lemma 5.2 in~\cite{surfacePDE}), such that
\begin{align}
    \frac{d}{dt}\int_{\Omega(t)}uv &= \int_{\Omega(t)}\big(v\partial_t u + u\partial_t v\big) + \int_{\Gamma(t)}\bm{\beta}\cdot\bm{n} u v\notag\\
    &= \int_{\Omega(t)}\big(v\partial_t u + u\partial_t v + \nabla\cdot(\bm{\beta} u v)\big)\notag\\
    &= \int_{\Omega(t)}\big(v\partial_t u + u\partial_t v + v\nabla\cdot(\bm{\beta}u) + u\bm{\beta}\cdot\nabla v\big).\label{eq:eee}
\end{align}
Combining~\eqref{eq:timebulk} with~\eqref{eq:eee}, integrating over $I_n$, integrating by parts the diffusion term, and adding stabilization terms results in the following scheme:

\begin{tcolorbox}
Given $u_h^- = u_h(t_{n-1}^-,\bm{x})$ find $u_h\in W_{h,k,m}^n$ such that
\begin{equation}
    A^n(u_h,v_h) + S^n(u_h,v_h) = L^n(v_h),\quad \text{for all }v_h\in W_{h,k,m}^n,\label{eq:conservative-scheme}
\end{equation}
where
\begin{align}
\begin{split}
    A^n(u, v) &= (u, v)_{\Omega(t_{n})} - \int_{I_n}(u,\partial_t v)_{\Omega(t)} - \int_{I_n}(u, \bm{\beta}\cdot \nabla v)_{\Omega(t)}\\
    &+ \int_{I_n}(D\nabla u, \nabla v)_{\Omega(t)},\label{eq:Ah-bulk-conservative}
\end{split}\\[1ex]
    S^n(u,v) &= \int_{I_n}s_h^n(t, u, v),\\
    L^n(v) &= (u_h^-, v)_{\Omega(t_{n-1})} + \int_{I_n}(f,v)_{\Omega(t)}.
\end{align}
\end{tcolorbox}
Equation~\eqref{eq:conservative-scheme} is solved one space-time slab at a time for $n=1,\dots, N$. The intitial condition, $u_h(t_{n-1}^-, \bm{x}) = \lim_{t \rightarrow t_{n-1}}u_h(t, \bm{x})$ is the solution from the previous space-time slab if $n>1$, and if $n=1$, $u_h(t_{n-1}^-, \bm{x}) = u_0(\bm{x})$. 
The stabilization terms $s_h^n$ are computed either as face-based ghost-penalty stabilization~\cite{ghost-penalty-stabilization, hansbo2016cut} or as patch-based ghost-penalty stabilization~\cite{preuss2018, lehrenfeld}. The face-based stabilization is defined as 
\begin{equation}
    s_h^n(t, u,v) = \sum_{i=1}^m\left(\sum_{F\in \F_h^n} \tau_F^i h^{2i-1} (\llbracket D_{\bm{n}_F}^i u\rrbracket_F, \llbracket D_{\bm{n}_F}^i v \rrbracket_F)_F \right),\label{eq:face-ghost-penalty}
\end{equation}
where $\tau_F^i>0$ for $i=1,\dots,m$, $D_{\bm{n}_F}^i$ denotes the $i$-th order derivative in the direction of the normal $\bm{n}_F$ to the face $F$, and $\llbracket \cdot \rrbracket_F$ denotes the jump operator across the face. In the standard ghost penalty stabilization which we refer to as full stabilization, the set of faces $\F_h^n$ is defined as 
\begin{align}
    \F_{h}^{n} &= \{F=K_1\cap K_2 : K_1,K_2\in\T_{h}^n,\ K_1\in\T_{h,\Gamma}^n\}.\label{eq:stabilization-faces-time}
\end{align}
The patch-based stabilization~\cite{preuss2018} is defined as
\begin{equation}
    s_h^n(u,v) = \sum_{F\in \F_h^n} \tau h^{-2}(\llbracket u\rrbracket_{\patch}, \llbracket v\rrbracket_{\patch})_{\patch},\label{eq:patch-ghost-penalty}
\end{equation}
where $\tau>0$, and for each $ F \in \F_h^n$, $\patch$ is the patch consisting of the two elements sharing face F and is defined by 
\begin{equation}
    \patch = K_1 \cup K_2,\ \text{where }K_1\cap K_2 = F \in \F_h^n.
\end{equation}
The patch jump is defined as 
\begin{equation}
    \llbracket u\rrbracket_{\patch}(\bm{x}) = u_1(\bm{x})-u_2(\bm{x}),
\end{equation} 
where the polynomial $u_i$ is defined as $u_i|_{K_i}=u|_{K_i}$ and outside $K_i$,  $u_i$ is the canonical extension of $u|_{K_i}$ so $u_i$ is well-defined for all $\bm{x}\in \patch$. The patch-based stabilization has an advantage when using high-order elements as it avoids the evaluation of derivatives.

In Section~\ref{section:stabilization}, we present an alternative set to $\F_h^n$ for the stabilization, based on constructing a macroelement partition of the active mesh.

By choosing $v_h=1$ in~\eqref{eq:conservative-scheme} we see that the scheme fulfills a discrete mass conservation equation:
\begin{equation}
    \int_{\Omega(t_n)}u_h - \int_{\Omega(t_{n-1})}u_h^- = \int_{t_{n-1}}^{t_n}\int_{\Omega(t)}f,
\end{equation}
and summing over all time instances $n=1,\dots,N$, we see that $u_h$ satisfies~\eqref{eq:conservation-law-non-zero}:
\begin{equation}
    \int_{\Omega(T)}u_h - \int_{\Omega(0)}u_0 = \int_0^T\int_{\Omega(t)}f.
\end{equation}

\subsection{Non-conservative scheme}

We also state the variational formulation used in e.g.~\cite{hansbo2016cut, lehrenfeld, trackers}, since we compare with this formulation in the numerical examples. On each time slab, we solve the following problem:
\begin{tcolorbox}
    Given $u_h^- = u_h(t_{n-1}^-,\bm{x})$ find $u_h\in W_{h,k,m}^n$ such that
\begin{equation}
    A^n(u_h,v_h) + S^n(u_h,v_h) = L^n(v_h),\quad \text{for all }v_h\in W_{h,k,m}^n,
    \label{eq:non-cons-eq}
\end{equation}
where
\begin{align}
\begin{split}
    A^n(u, v) &= (u, v)_{\Omega(t_{n-1})} + \int_{I_n}(\partial_t u, v)_{\Omega(t)} + \int_{I_n}(\bm{\beta}\cdot \nabla u, v)_{\Omega(t)}\\
    &+ \int_{I_n}(D\nabla u, \nabla v)_{\Omega(t)},\label{time:Ah-bulk}
\end{split}\\[1ex]
    S^n(u,v) &= \int_{I_n}s_h^n(t,u,v),\\
    L^n(v) &= (u_h^-, v)_{\Omega(t_{n-1})} + \int_{I_n}(f,v)_{\Omega(t)}.
\end{align}
\end{tcolorbox}

\section{Macroelement Stabilization}\label{section:stabilization}

Stabilization of the proposed cut finite element method is necessary since we want to have control of functions on the entire active mesh. In the next section, we will replace integrals over $I_n$ with quadrature rules. Without ghost-penalty stabilization, the scheme may lack control of functions in $W_{h,k,m}^n$ in the entire active mesh. Depending on how the interface cuts the background mesh, situations may occur where there are basis functions in the active finite element space that have arbitrarily small support in $\Omega(t_q)$ for all quadrature points $t_q$. Computationally, this leads to a system matrix with an unbounded condition number. 

To control the condition number, ghost-penalty stabilization using either the face-based technique defined in~\eqref{eq:face-ghost-penalty} or the patch-based technique defined in~\eqref{eq:patch-ghost-penalty} is used. Common for both techniques is that they stabilize elements connected to the faces in the set $\F_h^n$ defined in~\eqref{eq:stabilization-faces-time}. Stabilizing all faces (or face-patches) in this set is here referred to as full stabilization, see e.g.~\cite{hansbo2016cut, lehrenfeld, frachon2022cut}.

In~\cite{larson2021conservative}, a macroelement ghost-penalty stabilization is proposed for convection-diffusion equations in stationary domains. The main idea is to construct macroelements with a large intersection with the domain $\Omega$ and thus never stabilize between macro elements but only inside these elements. Here, we extend this idea to time-dependent domains. 

Let $Q_n = \{t_q^n\}_{q=1}^{N_t}$ denote a set of $N_t$ time instances in the time interval $I_n$. We classify an element $K$ in the active mesh $\T_{h}^n$ as large if 
\begin{equation}\label{eq:largeel}
    \delta \leq \frac{|K\cap \Omega(t_q)|}{|K|},  \quad \text{for all } t_q\in Q_n. 
\end{equation}
 Here $\delta\in (0,1]$. We specify in the numerical examples how we choose $\delta$. In the proposed method, we choose the set $Q_n$ to consist of the quadrature points associated with the quadrature rule we use to approximate integrals over $I_n$. See the next section. If an element is not large, it is classified as small. 

A macroelement partition $\M_{h}^n$ of $\Omega_{\T_{h}^n}$ is then constructed by the following algorithm.  

\begin{itemize}
    \item Let 
    \begin{equation}
        \T_L^n = \{K\in \T_h^n : |K\cap\Omega(t_q)|\geq \delta |K| \text{ for all }t_q\in Q_n\},
    \end{equation}
    be the set of large elements in the active mesh $\T_h^n$. 

    \item To each large element $K_L \in \T_L^n$ we associate a macroelement mesh
    \begin{align}
    \begin{split}
        \T_h^n(K_L) = \{&K \in \T_h^n : K = K_L \text{ or } K\in \T_h^n\backslash \T_L^n \text{ is face connected to }\\
        &K_L \text{ via a bounded number of internal faces}\}.
    \end{split}\label{eq:macroelement_mesh}
    \end{align}

    \item Each element $K\in \T_h^n$ belongs to precisely one macroelement mesh $\T_h^n(K_L)$.
    
    \item A macro element $M_L\in\M_{h}^n$ is defined as
    \begin{equation}
        M_L = \bigcup_{K\in \T_h^n(K_L)}K
    \end{equation}
    and the macroelement partition of $\Omega_{\T_h^n}$ is defined as
    \begin{equation}
        \M_{h}^n = \bigcup_{M\in\M_{h}^n}\{M\}.
    \end{equation}
\end{itemize}
Let $\F_{h}^n(M)$ denote the set of internal faces of the macroelement $M_L\in\M_{h}^n$. Hence, $\F_h^n(M)$ is empty if $\T_h^n(K_L)$ only consists of the large element $K_L$. The faces to which stabilization terms are added in the method are given by  
\begin{align}
    \F_{h}^{\M,n} = \bigcup_{M\in\M_{h}^n}\F_{h}^n(M).
\end{align}
In the stabilization of the method using~\eqref{eq:face-ghost-penalty} or~\eqref{eq:patch-ghost-penalty}, $\F_h^n$ is replaced by $\F_h^{\M,n}$. Thus, only the internal faces of the macroelements are considered in the stabilization process. 

\begin{figure}[H]
    \centering
    \includegraphics[width=.475\textwidth]{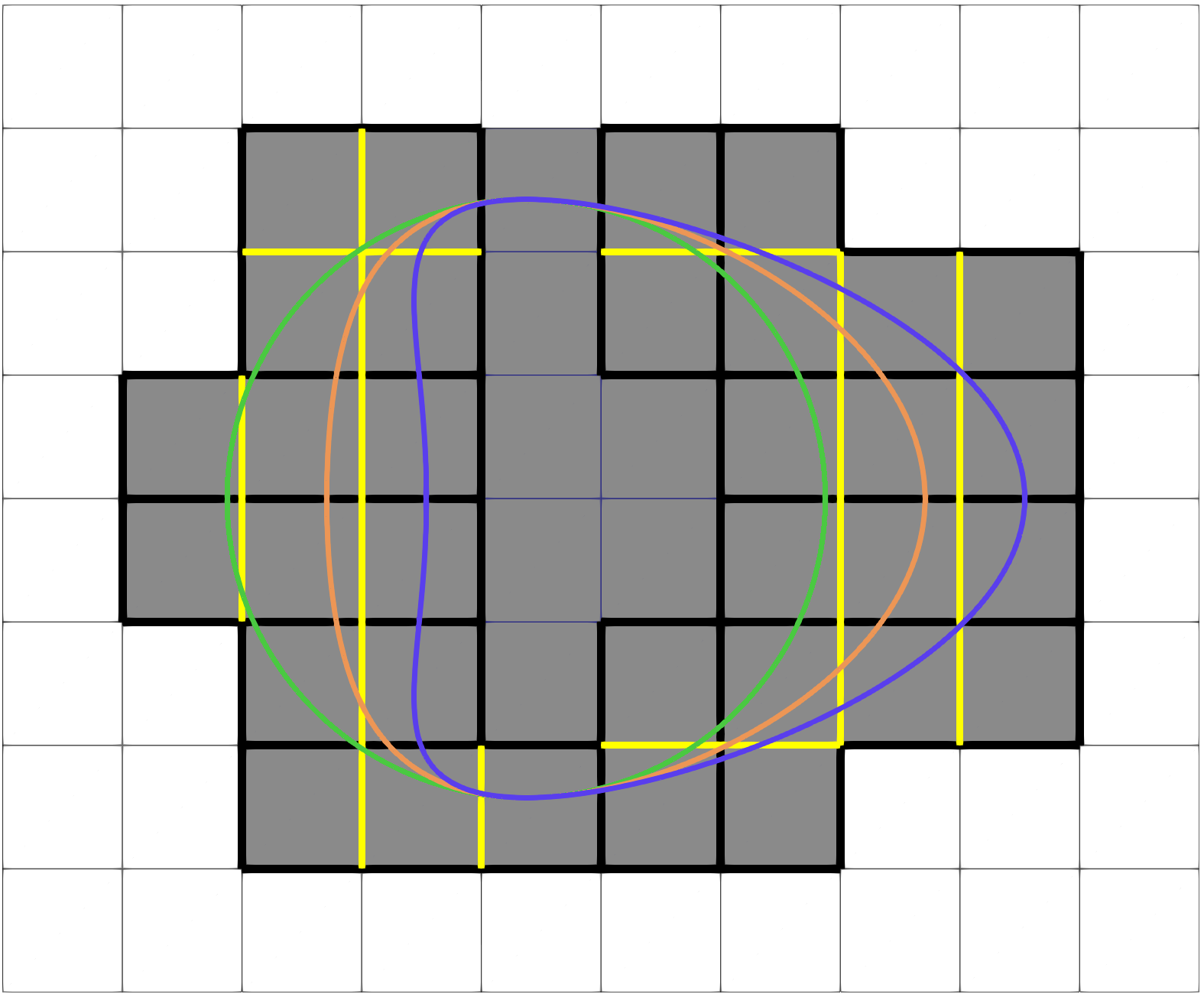}%
    \hfill
    \includegraphics[width=.475\textwidth]{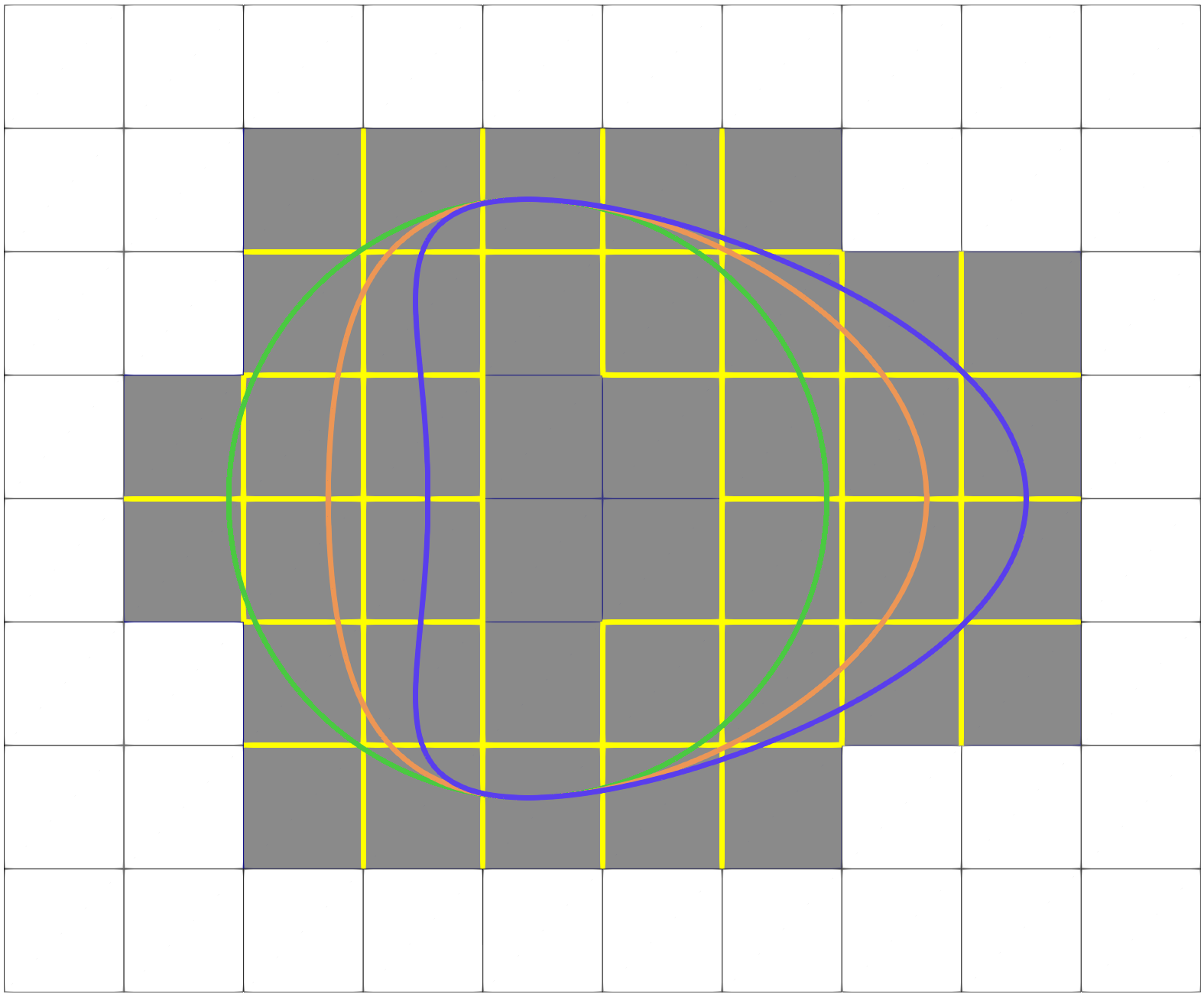}\caption{Macroelement stabilization (left) versus full stabilization (right). The surface is here visualized at three time instances, $Q_n = \{t_{n-1}, (t_{n-1}+t_n)/2, t_n\} \subset I_n$ which are used to construct the macroelement partition. Black faces mark external borders of macroelements $M\in\M_{h}^n$, and the yellow faces mark the internal faces $F\in \F_{h}^{\M,n}$ in the left panel and $F\in \F_h^{n}$ in the right panel. The parameter $\delta$ is chosen as $\delta = 0.3$.}\label{fig:fullstab_vs_macro}
\end{figure}

For an illustration of the macroelement partition and its corresponding stabilized faces compared to full stabilization, see Figure~\ref{fig:fullstab_vs_macro}.  It is noticeable that using a macroelement partition reduces the number of stabilized faces significantly, leading to greater sparsity of the system matrix. (Compare the number of yellow faces in the left panel of Figure~\ref{fig:fullstab_vs_macro} with the right panel.) Note that the construction of the macroelement mesh~\eqref{eq:macroelement_mesh} is not unique, since there are several ways a small element can be connected to a large element.

\section{Numerical Quadrature}\label{section:numerical_quadrature}
We approximate the integrals in the variational formulations in Section~\ref{section:convection-diffusion} by first using a quadrature rule for the integral in time. In this way we never explicitly construct the space-time domain in $\R^{d+1}$. In this section, we discuss how these integrals are approximated using quadrature rules that allow for high-order approximations.

\subsection{Quadrature in time}\label{subsection:numerical-intergration}
In this work, integrals in time are approximated by the Gauss-Lobatto quadrature rule $\{t_q^n, \weight_q\}_{q=1}^{N_t}$, where $\{t_q^n\}_{q=1}^{N_t}$ are evenly spaced points in the time interval $I_n$ such that $t_1^n = t_{n-1}$ and $t_{N_t}^n = t_n$, and $\{\weight_q\}_{q=1}^{N_t}$ are the corresponding Gauss-Lobatto weights. For a sufficiently smooth $f$, we have
\begin{equation}
    \int_{I_n} f(t,\bm{x}) dt = \sum_{q=1}^{N_t}\weight_q f(t_q^n,\bm{x})+ \mathcal{O}((\Delta t_n)^{2 N_t}).
\end{equation}

The proposed conservative scheme then becomes: Given $u_h^- = u_h(t_{n-1}^-,\bm{x})$ find $u_h\in W_{h,k,m}^n$ such that
\begin{equation}
    A_h^n(u_h,v_h) + S_h^n(u_h,v_h) = L_h^n(v_h),\quad \text{for all }v_h\in W_{h,k,m}^n,\label{eq:conservative-scheme-quadrature}
\end{equation}
where
\begin{align}
\begin{split}
    A_h^n(u, v) &= (u, v)_{\Omega(t_{n})} - \sum_{q=1}^{N_t}\weight_q (u,\partial_t v)_{\Omega(t_q)} - \sum_{q=1}^{N_t}\weight_q (u, \bm{\beta}\cdot \nabla v)_{\Omega(t_q)}\\
    &+ \sum_{q=1}^{N_t}\weight_q(D\nabla u, \nabla v)_{\Omega(t_q)},\label{eq:Ah-bulk-conservative-quadrature}
\end{split}\\[1ex]
    S_h^n(u,v) &= \sum_{q=1}^{N_t}\weight_q s_h^n(t_q, u, v),\\
    L_{h}^n(v) &= (u_h^-, v)_{\Omega(t_{n-1})} + \sum_{q=1}^{N_t}\weight_q(f,v)_{\Omega(t_q)}.
\end{align}
The number of quadrature points $N_t$ is chosen to be the same for all intervals $I_n$. In the numerical experiments, we study the number of quadrature points $N_t$ that have to be used for an optimal rate of convergence of the method. 

\subsection{Quadrature in space}\label{subsection:quadrature_cut}

A challenge in CutFEM is to accurately compute integrals on elements that are cut by the boundary.
In this work, a level set function $\phi: I\times \Omega_0 \rightarrow \R$ is used to define the domain and its boundary in any given time instance by
\begin{align}
    \Omega(t) &= \{\bm{x}\in \Omega_0 : \phi(t,\bm{x}) < 0\},\\
    \Gamma(t) &= \{\bm{x}\in \Omega_0 : \phi(t, \bm{x}) = 0\}.
\end{align}
We use the algorithm first presented in~\cite{saye} to approximate space integrals. It gives tensor product quadrature rules using $N_s$ nodes in each space dimension resulting in a rule $\{\bm{x}_q, \alpha_q\}_{q=1}^{N_s^d}$ for integrals on elements cut by the bulk domain, $\Omega(t)\cap K$, and a rule $\{\bm{s}_q, \gamma_q\}_{q=1}^{N_s^{d-1}}$ for integrals on elements cut by the surface $\Gamma(t)\cap K$. The quadrature rule $\sum_{q=1}^{N_s^d}\alpha_q g(\bm{x}_q)$ approximates integrals $\int_{\Omega(t)\cap K} g(\bm{x}) d\bm{x}$ and $\sum_{q=1}^{N_s^{d-1}}\gamma_q g(\bm{s}_q)$ approximates $\int_{\Gamma(t)\cap K} g(\bm{s}) d\bm{s}$ with order of accuracy of approximately $2 N_s$, inheriting the accuracy of standard Gaussian quadrature rules~\cite{saye}.

\begin{figure}[H]
    \centering
    \includegraphics[width=.35\linewidth]{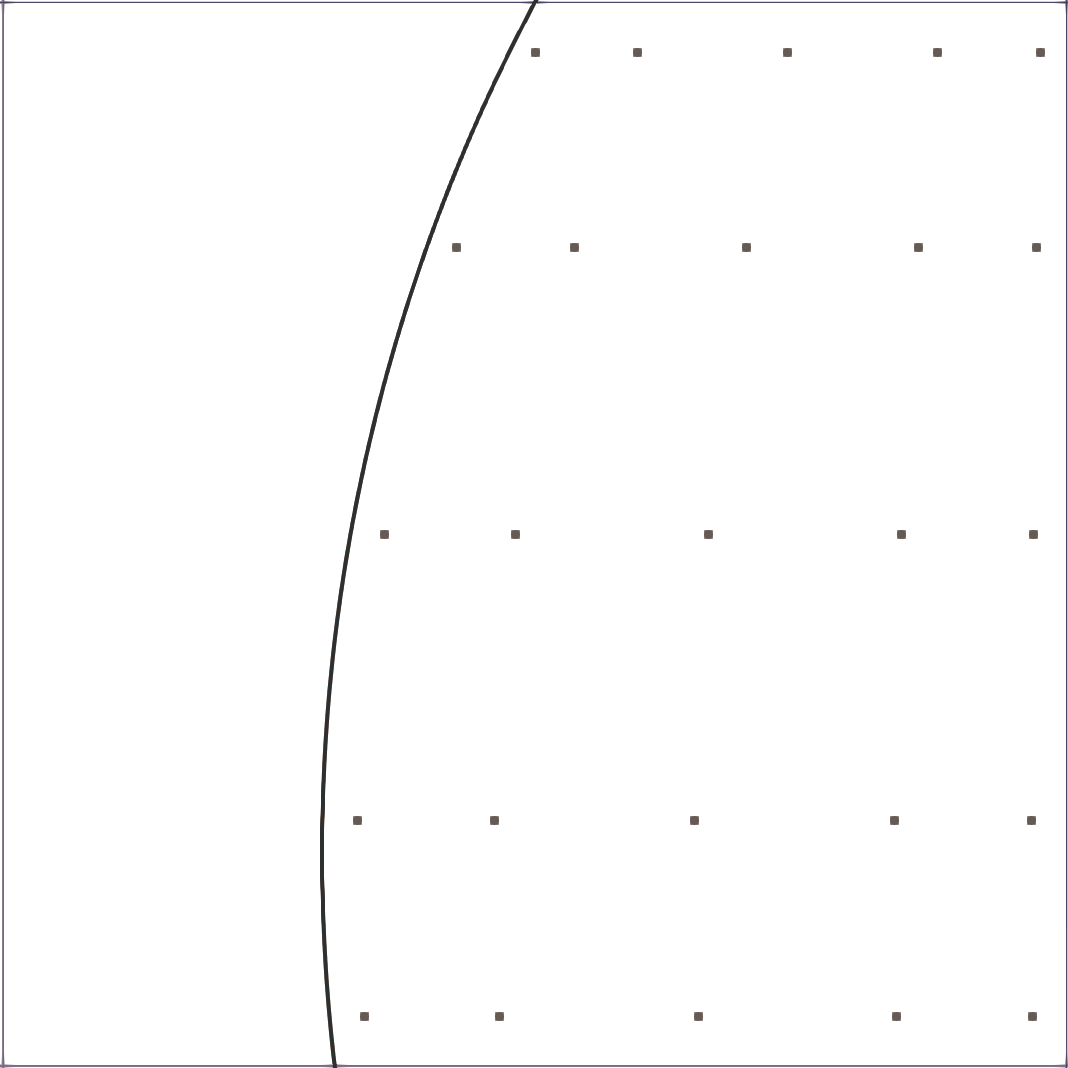}
    \caption{Quadrature nodes produced by the algorithm presented in~\cite{saye} for approximating an integral on a cut element $\Omega(t)\cap K$. The nodes are distributed as a tensor product between two one-dimensional quadrature rules with five quadrature nodes. This results in an order of accuracy of $10$.}\label{fig:saye}
\end{figure}

The algorithm is based on transforming the zero contour of the level set function into the graph of a height function in the coordinate system of the two spatial directions. The existence of such a height function is in most cases guaranteed by the implicit function theorem, and in the cases where its conditions are not fulfilled, the element in question is split into several elements. The recursive algorithm leads to a tensor product Gaussian quadrature rule of high order, and quadrature points for a $10$th order accurate rule are shown in Figure~\ref{fig:saye}. For a more elaborate algorithm description, see~\cite{saye}.

\section{Numerical Examples}\label{section:numerical-examples}

In this section, we present two numerical examples and examine the methods presented in the previous sections. In particular, we compare the performance of the non-conservative scheme with the conservative scheme. All results presented in this section are with patch-based stabilization. However, we also tested the standard face-based ghost penalty stabilization and the results are similar. The code that was used to run the numerical experiments is an in-house code written in C++, based on the open-source FreeFEM++ library~\cite{freefem}. The numerical quadrature discussed in Section~\ref{subsection:quadrature_cut} is based on the open-source library Algoim, developed by R. Saye~\cite{algoim}, which was integrated into our in-house code. For reproducibility, an open-source version of this code is available~\cite{github}.

We measure the error using the following norms.
\begin{align}
    \|u-u_h\|_{L^2(\Omega(T))} &= \int_{\Omega(T)}|u(T,x)-u_h(T,x)|^2 dx,\label{eq:norm1}\\
    \|u-u_h\|_{L^2(L^2(\Omega(t)), 0, T)}^2 &= \int_{I_n}\|u-u_h\|_{L^2(\Omega(t))}^2 dt.\label{eq:norm2}
\end{align}
We measure the conservation error as
\begin{equation}
    e_c(T) = \Bigg|\int_{\Omega(T)}u_h(T,\bm{x}) - \int_{\Omega(0)}u_0(\bm{x}) - \sum_{n=1}^N \int_{I_n} \int_{\Omega(t)} f(t,\bm{x}) \Bigg|.
    \label{eq:cons-global-error-bulk}
\end{equation}

\subsection{A moving circle}\label{subsection:ex1}
\begin{figure}[H]
    \centering
    \includegraphics[width=1.\textwidth]{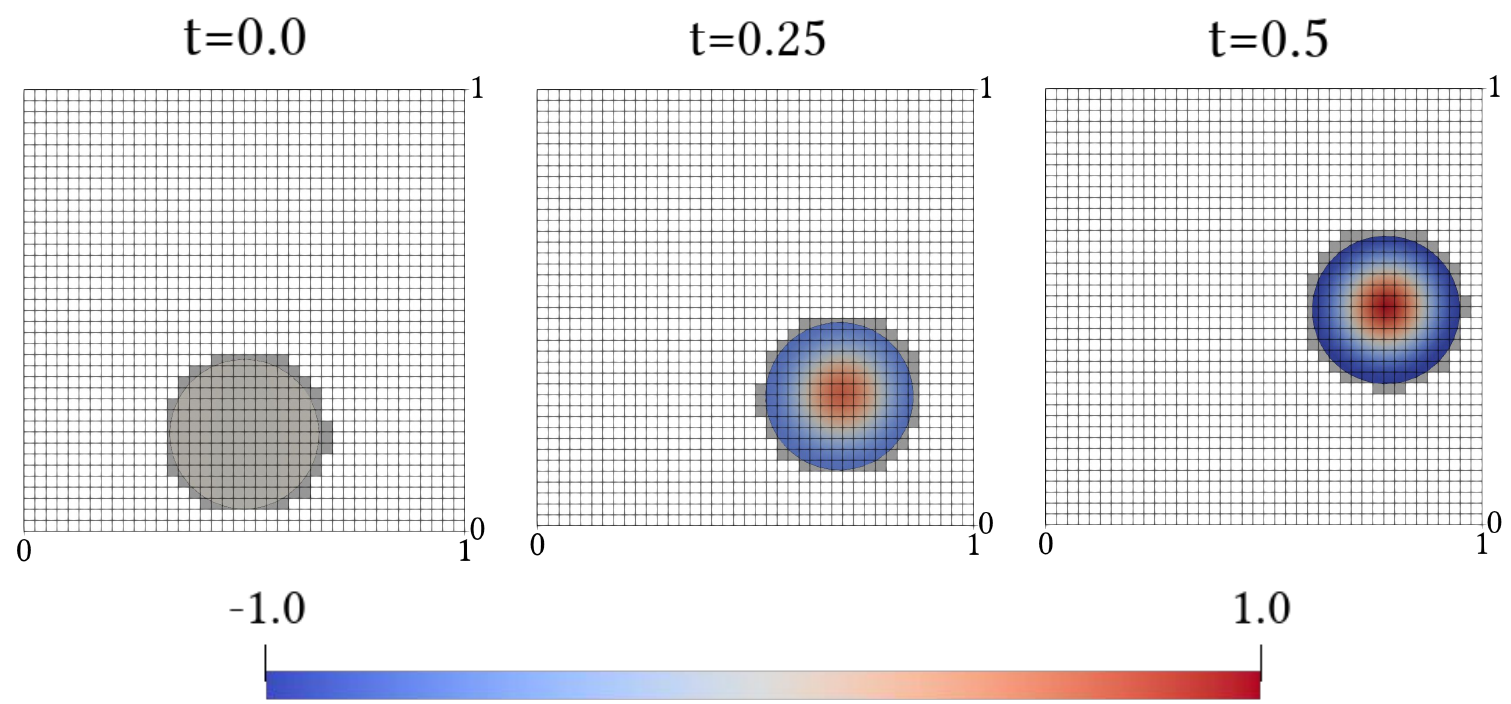}
    \caption{Example~\ref{subsection:ex1}: Numerical solution using $h=0.025$. The active mesh is illustrated in grey.}\label{fig:ex1-sol}
\end{figure}

The boundary evolves with the velocity field $\bm{\beta}=(\pi(0.5-y), \pi(x-0.5))$ and the distance function is explicitly given by 
\begin{equation}
    \phi(t,x,y) = r(t,x,y)^2 - r_0^2, \label{eq:level_set}  
\end{equation}
where $r(t,x,y)=\sqrt{(x-x_c(t))^2+(y-y_c(t))^2}$, $x_c(t) = 0.5 + 0.28\sin(\pi t)$ and $y_c(t) = 0.5 - 0.28\cos(\pi t)$. In this example, $r_0=0.17$. The diffusion coefficient is chosen as $D=1$, and the source function $f$ is computed such that the exact solution $u$ of~\eqref{eq:timebulk}--~\eqref{eq:neumann} is
\begin{equation}
    u(t,x,y) = \cos(\pi r(t,x,y)/r_0)\sin(\pi t),\quad (t,x,y)\in I\times \Omega(t).\label{eq:exact_solution}
\end{equation}
The numerical solution is illustrated on a uniform mesh in Figure~\ref{fig:ex1-sol}. Unless stated otherwise, the time step size is set to $\Delta t_n = h/3$.

\begin{figure}[H]
    \centering
    \includegraphics[width=.8\linewidth]{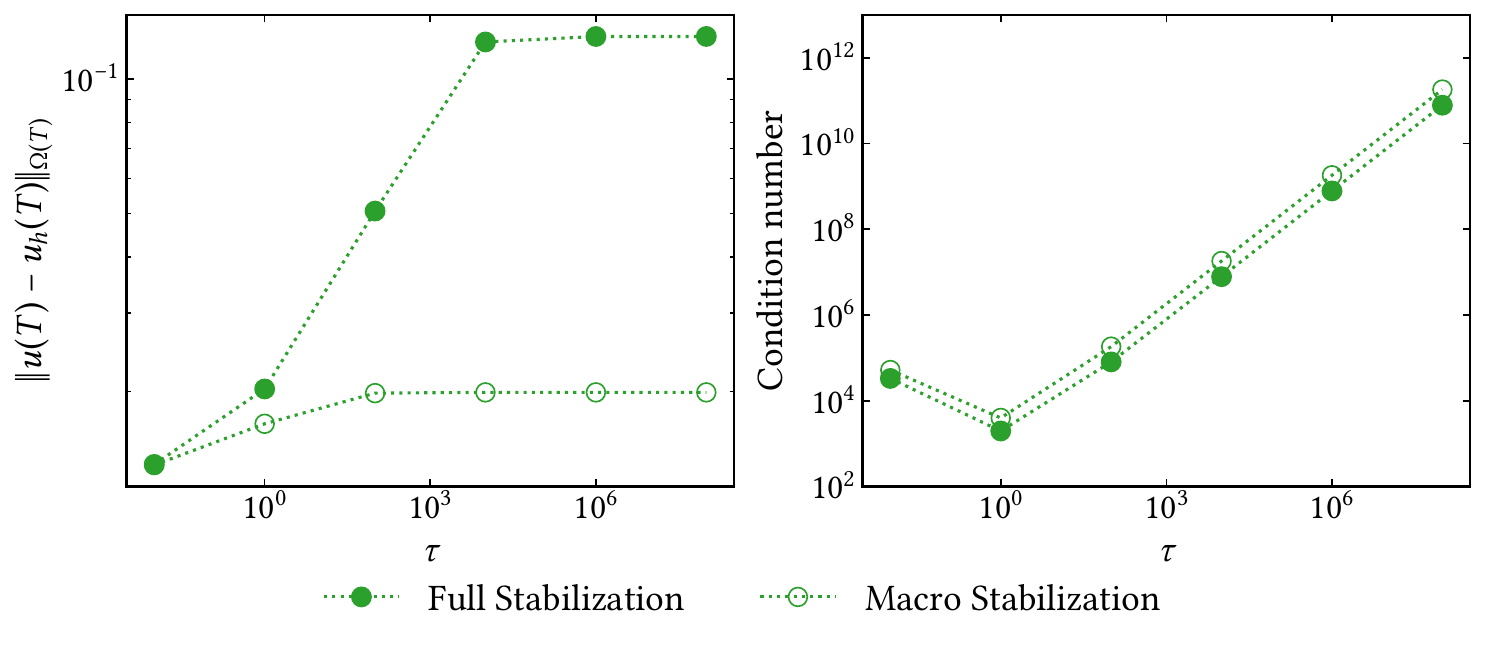}
    \caption{Example~\ref{subsection:ex1}. $L^2$ error and condition number versus the stabilization parameter $\tau$. The conservative scheme is used with linear elements in both space and time ($m=k=1$), with $h=0.05$, $T=0.5$, and $\delta=0.5$. The full and macroelement stabilization techniques are compared.}\label{fig:L2_vs_tau_example1}
\end{figure}

We study the macroelement stabilization with respect to the stabilization parameter $\tau$ and the macroelement parameter $\delta$ and compare it with full stabilization. In Figure~\ref{fig:L2_vs_tau_example1}, the $L^2$ error (measured in the last time instance), and the condition number of the corresponding matrix, are shown versus $\tau$. Linear elements are used in both space and time, $m=k=1$.  It is clear that with the macroelement stabilization, the $L^2$ error is more or less invariant as the stabilization parameter increases, meanwhile with full stabilization, the error grows significantly. The condition number is slightly larger when using macroelement stabilization but follows the same trend as with full stabilization. 

\begin{figure}[H]
    \centering
    \includegraphics[width=.95\linewidth]{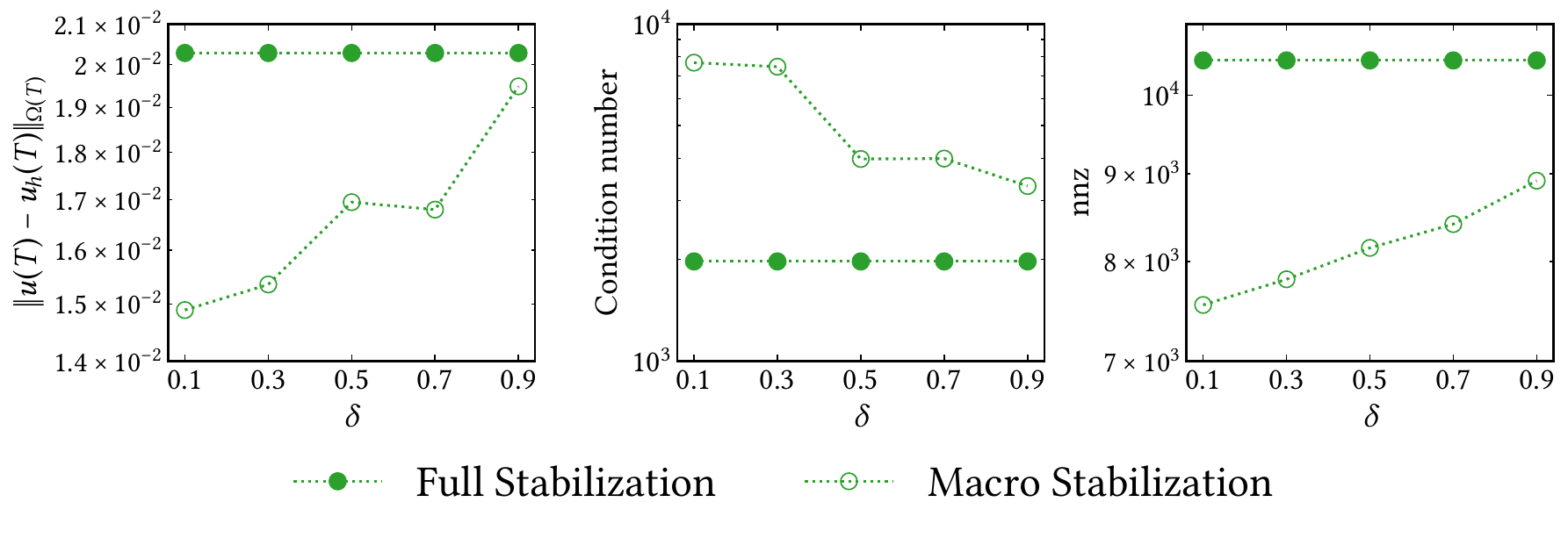}
    \caption{Example~\ref{subsection:ex1}. $L^2$ error, condition number, and the number of non-zero elements versus the macroelement parameter $\delta$. The conservative scheme is used using linear elements in both space and time, i.e., $m=k=1$, with $h=0.05$, $T=0.5$, and $\tau=1$. The full and macroelement stabilization techniques are compared.}\label{fig:example1_vs_delta}
\end{figure}

In Figure~\ref{fig:example1_vs_delta}, the $L^2$ error, condition number, and the number of non-zero elements (nnz) in the matrix are plotted against the $\delta$-parameter. As $\delta$ increases and we stabilize more, the $L^2$ error increases while the condition number decreases, the number of non-zero elements increases, and the macroelement stabilization approaches full stabilization.

\begin{figure}
    \centering
    \includegraphics[width=1\linewidth]{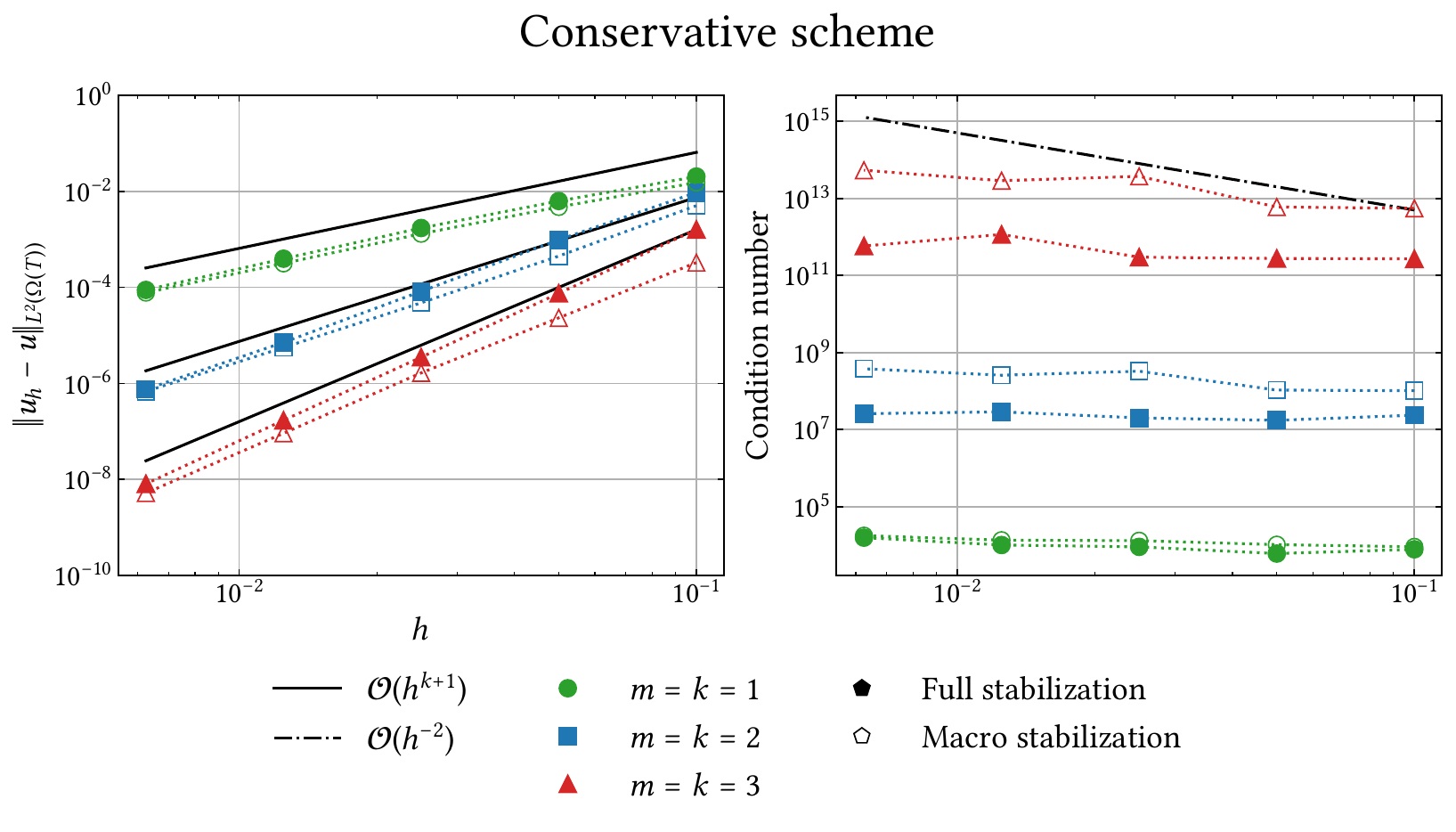}
    \caption{Example~\ref{subsection:ex1}: $L^2$ error and the condition number versus $h$ for the conservative scheme with $T=0.1$. For $m=k=1$, we use $N_t=3$, for $m=k=2$ we use $N_t=5$, and for $m=k=3$ we use $N_t=9$.}\label{fig:example1_convergence}
\end{figure}

\begin{figure}
    \centering
    \includegraphics[width=1\linewidth]{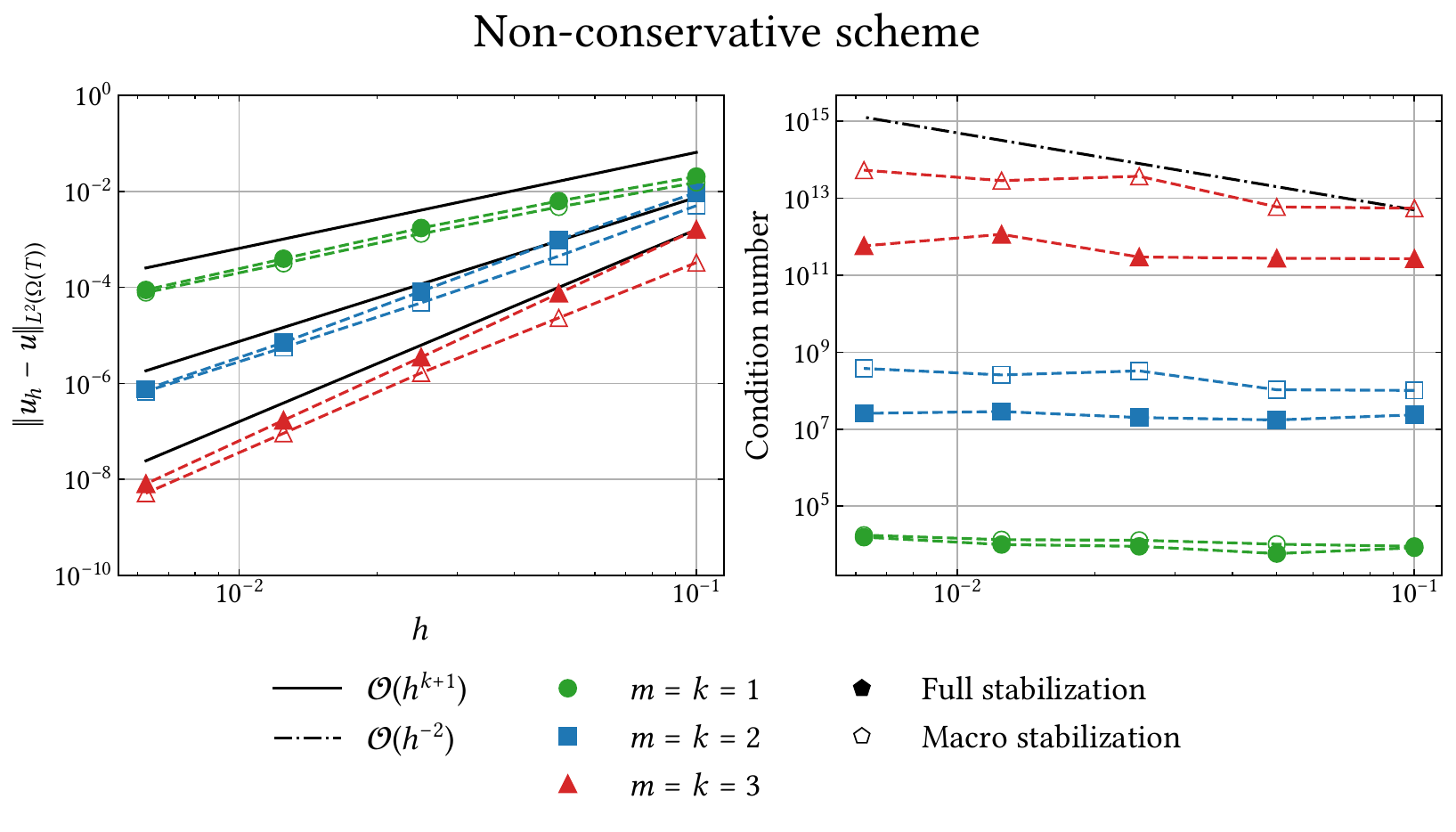}
    \caption{Example~\ref{subsection:ex1}: $L^2$ error and the condition number versus $h$ for the non-conservative scheme with $T=0.1$. For $m=k=1$, we use $N_t=3$, for $m=k=2$ we use $N_t=5$, and for $m=k=3$ we use $N_t=9$.}\label{fig:example1_convergence_classical}
\end{figure}

Based on the results of Figure~\ref{fig:L2_vs_tau_example1} and Figure~\ref{fig:example1_vs_delta}, we choose $\tau=1$ and $\delta = 0.5$ for the rest of the numerical experiments in this example unless stated otherwise. We emphasize that smaller $\delta$ could also be used. In Figure~\ref{fig:example1_convergence} and Figure~\ref{fig:example1_convergence_classical} we investigate the convergence of the method with respect to the mesh size, comparing the conservative and non-conservative schemes. The results show that the conservative and non-conservative schemes both converge optimally, i.e.\ the order of convergence is $k+1$ for the $L^2$ error using polynomial orders $m=k$. The condition number is stable using both macro and full stabilization, growing slower with respect to decreasing $h$ than the theoretically expected $\mathcal{O}(h^{-2})$.

\begin{figure}
    \centering
    \includegraphics[width=1\linewidth]{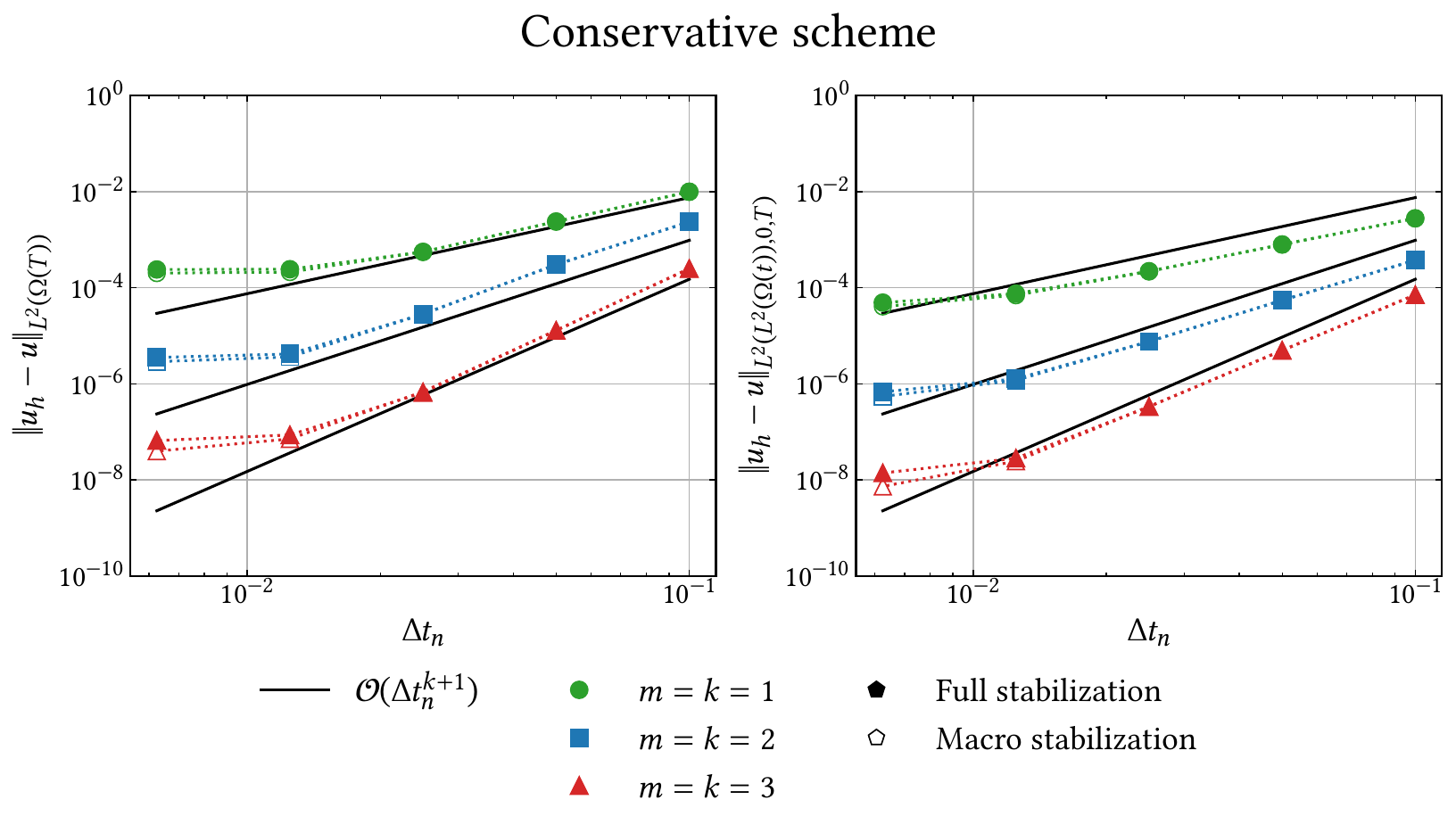}
    \caption{Example~\ref{subsection:ex1}: $L^2$ error versus $\Delta t_n$ for the conservative scheme with $T=0.1$ and $h=0.001$.}\label{fig:example1_convergence_vs_dt}
\end{figure}

\begin{figure}
    \centering
    \includegraphics[width=1\linewidth]{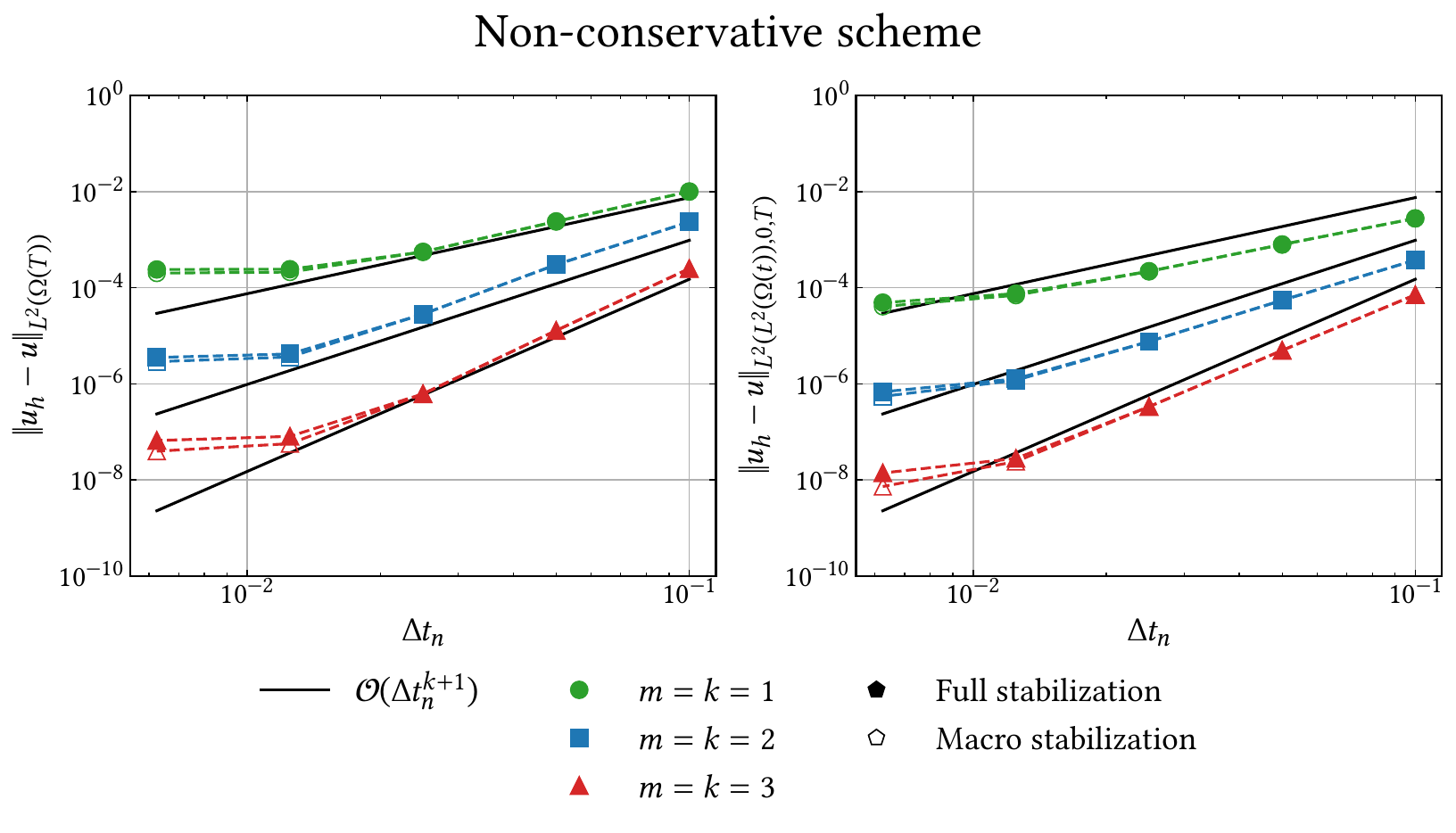}
    \caption{Example~\ref{subsection:ex1}: $L^2$ error versus $\Delta t_n$ for the non-conservative scheme with $T=0.1$ and $h=0.001$.}\label{fig:example1_convergence_vs_dt_classical}
\end{figure}

In Figure~\ref{fig:example1_convergence_vs_dt} and Figure~\ref{fig:example1_convergence_vs_dt_classical}, the convergence with respect to the time step size, measured in the two norms~\eqref{eq:norm1} and~\eqref{eq:norm2}, is illustrated. The mesh size is kept fixed at $h=0.001$. Both the conservative and the non-conservative schemes show convergence of order $k+1$ when polynomials of order $k$ are used in time.

\begin{figure}[H]
    \centering
    \includegraphics[width=1\linewidth]{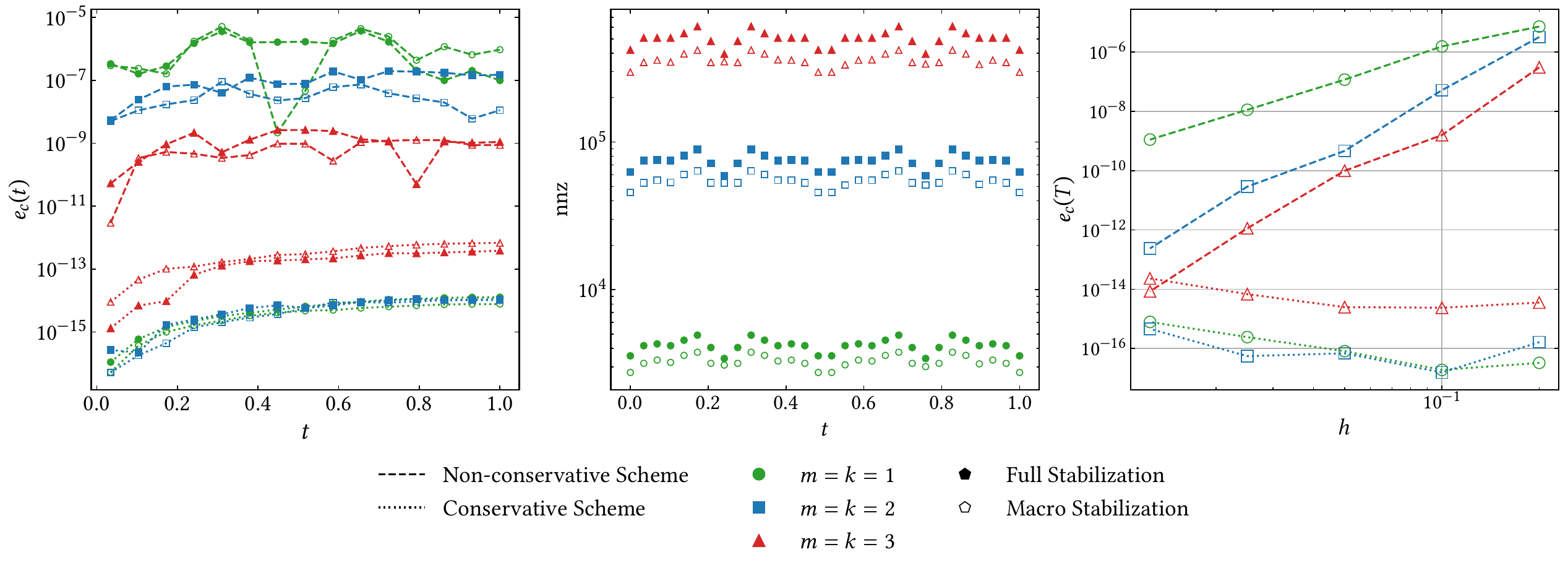}
    \caption{Example~\ref{subsection:ex1}: Conservation error $e_c(t)$ (left) and number of non-zero entries in the system matrix (middle) versus time $t$ with $h=0.1$. Comparison of the conservation error $e_c(T)$ for the conservative and non-conservative methods using macroelement stabilization (right). The error is plotted against $h$ with $T=0.1$.}\label{fig:cons-err-example1}
\end{figure}

The global conservation error $e_c(t)$ with respect to time is illustrated in Figure~\ref{fig:cons-err-example1}, showing that the conservation error from the conservative scheme is of the order of machine epsilon for any polynomial order, while this is not the case for the non-conservative solution. The figure also shows that the number of non-zero entries in the system matrix is reduced for all polynomial orders when using macroelement stabilization. Furthermore, the conservation error in time $T$ is shown versus the mesh size. From the figure, it is clear that the conservation error of the non-conservative scheme converges (as expected) with decreasing mesh size. However, one has conservation at machine precision order for any mesh size, and in particular, on very coarse meshes, with the conservative method.

\subsection{The kite}\label{subsection:ex2}

\begin{figure}[H]
    \centering
    \includegraphics[width=1\linewidth]{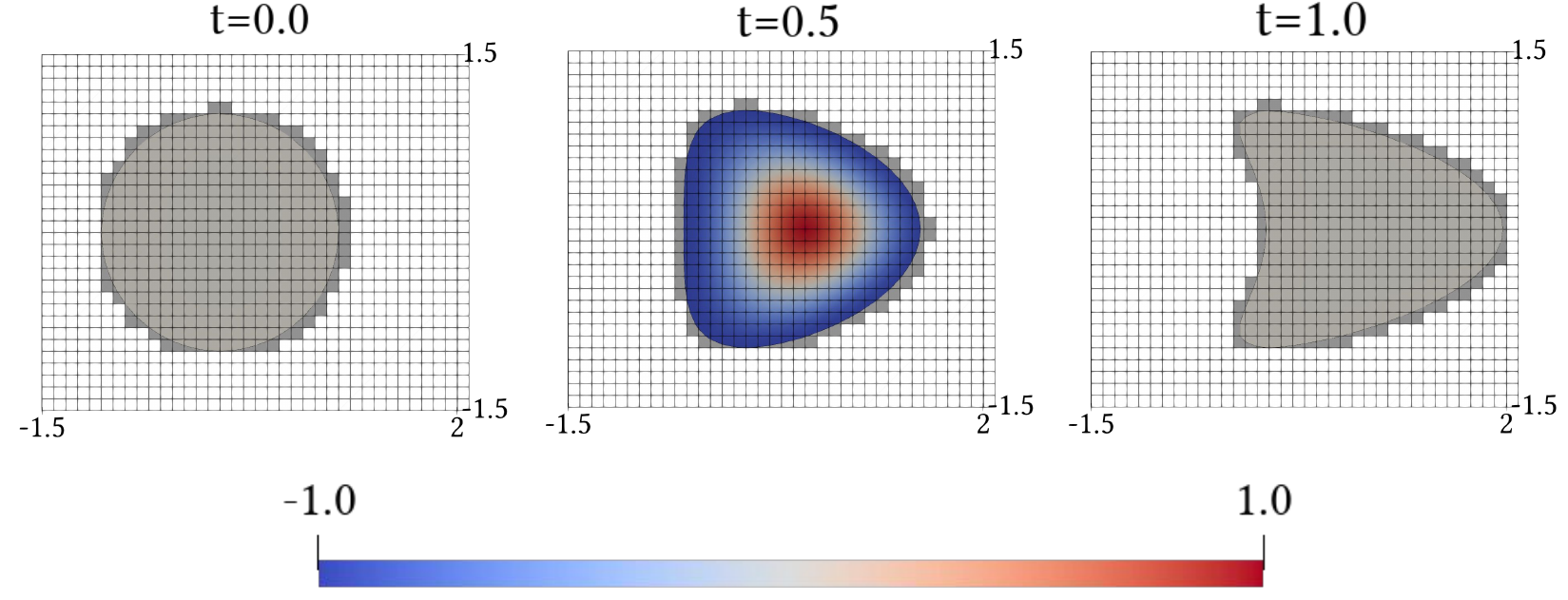}
    \caption{Example~\ref{subsection:ex2}: Numerical solution using $h=0.05$.}\label{fig:ex2-sol}
\end{figure}

This example is taken from Section 6.1 in~\cite{lehrenfeld}. The velocity field is given by $\bm{\beta} = (1-y^2, 0)$ and the boundary and the solution are as in equation~\eqref{eq:level_set} and~\eqref{eq:exact_solution} respectively, but with $x_c(t)=(1-y^2)t$ and $y_c(t)=0$. The diffusion coefficient is $D=1$. We set the stabilization parameter to $\tau= 0.1$, the macroelement parameter to $\delta=0.3$, and $\Delta t_n = 5h/18$ for the numerical experiments unless indicated otherwise. The numerical solution is illustrated in Figure~\ref{fig:ex2-sol}. 

\begin{figure}[H]
    \centering
    \includegraphics[width=1\textwidth]{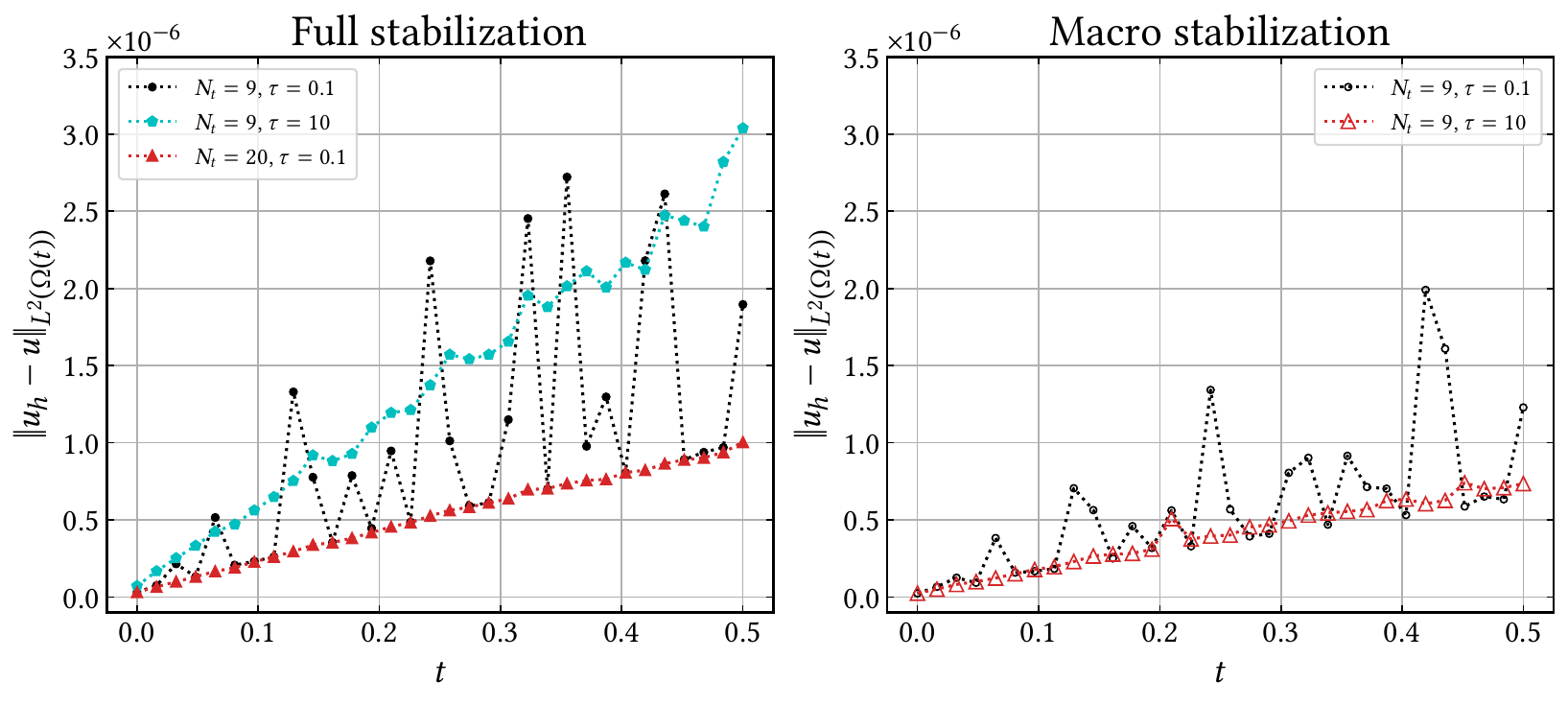}
    \caption{Example~\ref{subsection:ex2}: Increased stability of the conservative scheme with $m=k=3$ by using either a larger number of quadrature points in time using full stabilization (left) or larger stabilization constant using macroelement stabilization (right). The mesh size is $h = 0.05625$.}\label{fig:example2_error_vs_t}
\end{figure}

We studied how the stability of the conservative method was affected by the number of quadrature points in time and the stabilization constant, and found that the stability is more sensitive for high-order polynomials. In Figure~\ref{fig:example2_error_vs_t}, we illustrate the result for cubic polynomials in both space and time. We see that increasing the number of quadrature points in time or increasing the stabilization constant makes the $L^2$ error less sensitive to how elements in the background mesh are cut. However, increasing the stabilization constant increases the $L^2$ error when full stabilization is used, while this is not the case when using the macroelement stabilization we propose. Based on this observation, we choose $N_t=9$ and $\tau=10$ for the conservative method when using macroelement stabilization for $m=k=3$. For the method using full stabilization, we use $N_t=20$ and $\tau=0.1$. The macroelement stabilization is the preferable option, since it is more efficient, increases the sparsity of the system matrix, and the $L^2$ error is not increased by using a larger stabilization constant. 

\begin{figure}
    \centering
    \includegraphics[width=1\linewidth]{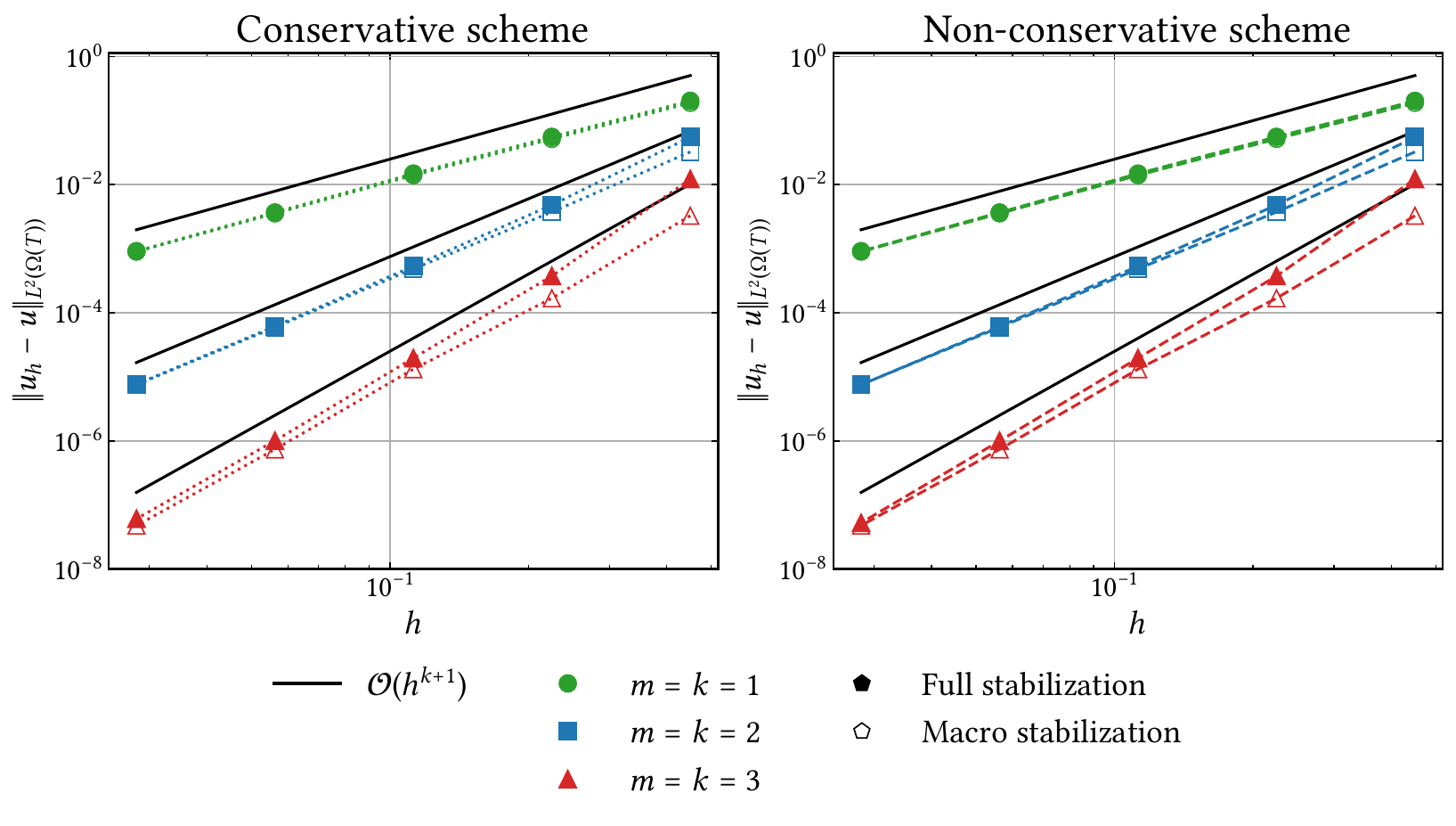}
    \caption{Example~\ref{subsection:ex2}: $L^2$ error versus $h$ at $T=0.5$.}\label{fig:example2_convergence}
\end{figure}

\begin{figure}
    \centering
    \includegraphics[width=.85\linewidth]{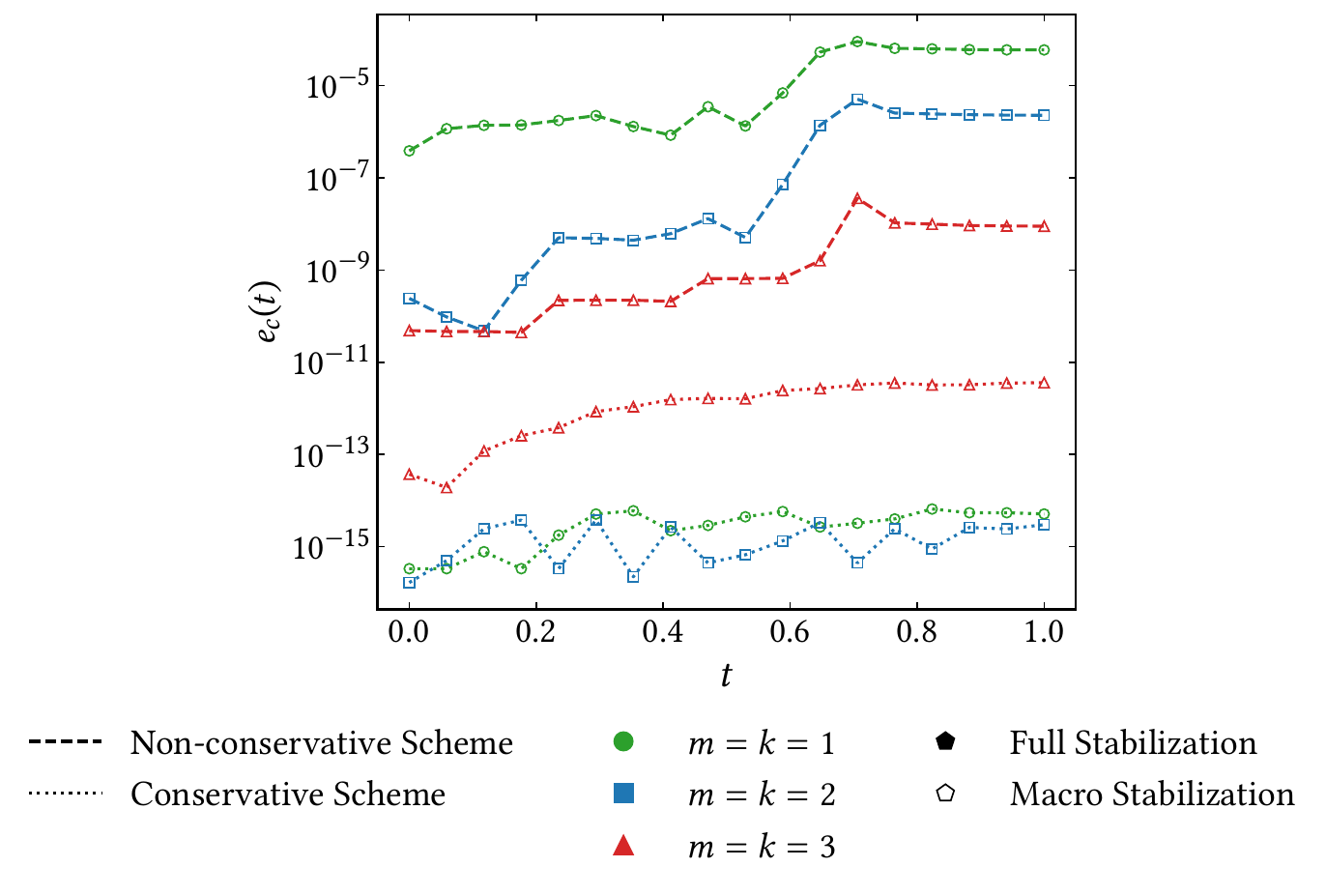}
    \caption{Example~\ref{subsection:ex2}: Conservation error $e_c(t)$ versus time $t$ with $h=0.2$ and $T=1$. The conservative and non-conservative schemes are shown using macroelement stabilization with $\delta = 0.3$.}\label{fig:cons-err-example2}
\end{figure}

The $L^2$ error versus mesh size is shown for both the conservative and the non-conservative scheme in Figure~\ref{fig:example2_convergence}. The conservation error $e_c(t)$ is illustrated in Figure~\ref{fig:cons-err-example2}. Both methods show optimal convergence, however, the conservation error of the conservative scheme is of the order of machine epsilon while the non-conservative scheme has a significantly larger conservation error. 

\section{Extension to Coupled Bulk-Surface Problems}\label{section:coupled_bulk_surface}

In this section, we extend the ideas presented in the previous sections to coupled bulk-surface convection-diffusion equations modeling the transport of soluble surfactants. We consider the following problem:

Find $u_B: I\times \Omega(t)\rightarrow \R$ and $u_S: I\times \Gamma(t)\rightarrow \R$ such that
\begin{alignat}{2}
    \partial_t u_B + \bm{\beta}\cdot\nabla u_B - \nabla\cdot(D\nabla u_B) &= f_B, \quad &&\text{in } I\times \Omega(t),\label{eq:bulk-coupled}\\
    \partial_t u_S + \bm{\beta}\cdot\nabla u_S + (\nabla_{\Gamma}\cdot\bm{\beta})u_S - \nabla_{\Gamma}\cdot(D_{\Gamma}\nabla_{\Gamma} u_S) &= f_S + f_C, \quad &&\text{on } I\times \Gamma(t),\label{eq:surface-coupled}\\
    -\bm{n}\cdot(D\nabla u_B) &= f_{C}, \quad &&\text{on } I\times \Gamma(t)\label{eq:coupling},
\end{alignat}
where $D$ and $D_{\Gamma}$ are the bulk and surface diffusion coefficients respectively, $f_B$ and $f_S$ are source terms, and initial conditions are assumed to be known. The operator $\nabla_{\Gamma}$ denotes the tangential gradient, see e.g.~\cite{sara_surface, burman2015}. The coupling term $f_C$ describes the exchange of surfactants between the bulk and the surface, and there are several models for this term, such as Henry, Langmuir, or Frumkin isotherms~\cite{hansbo2016cut, PROSSER20011}. In this work, we consider the Langmuir isotherm, which is a non-linear coupling term given by
\begin{equation}
    f_C = b_B u_B - b_S u_S - b_{BS} u_B u_S,
\end{equation}
where $b_B$, $b_S$, and $b_{BS}$ are physical constants related to the adsorption and desorption of the surfactant. The case of $b_{BS}=0$ gives the linear Henry coupling model, which is treated in e.g.~\cite{heimann2024higher}. In this work, we consider for brevity $b_B=b_S=b_{BS}=1$.

The finite element space used for this problem is given by
\begin{equation}
    W_{h,k,m}^n = \left(P_k(I_n)\otimes V_{h,m}|_{\Omega_{\T_h^n}}\right) \times \left(P_k(I_n)\otimes V_{h,m}|_{\Omega_{\T_{h,\Gamma}^n}}\right)
\end{equation}
and the discretizations are given as follows.

\subsection{Conservative scheme}

Let $u_h=(u_{B,h}, u_{S,h}) \in W_{h,k,m}^n$ denote the numerical solution. The method reads: Given $u_h^-=u_h(t_{n-1}^-, \bm{x})$ find $u_h\in W_{h,k,m}^n$ such that
\begin{equation}
    A_h^n(u_h,v_h)+S_h^n(u_h,v_h) = L_n^n(v_h),\quad \text{for all }v_h\in W_{h,k,m}^n,
\end{equation}
where
\begin{align}
    A_h^n(u, v) &= (u_B, v_B)_{\Omega(t_{n})} + (u_S,v_S)_{\Gamma(t_n)} - \int_{I_n}(u_B,\partial_t v_B)_{\Omega(t)}\nonumber\\ 
    &- \int_{I_n}(u_S, \partial_t v_S)_{\Gamma(t)} - \int_{I_n}(u_B, \bm{\beta}\cdot \nabla v_B)_{\Omega(t)} - \int_{I_n}(u_S, \bm{\beta}\cdot\nabla v_S)_{\Gamma(t)}\nonumber\\ 
    &+ \int_{I_n}(D\nabla u_B, \nabla v_B)_{\Omega(t)} + \int_{I_n}(D_{\Gamma}\nabla_{\Gamma} u_S, \nabla_{\Gamma} v_S)_{\Gamma(t)}\label{eq:A-coupled-c}\\
    &+ \int_{I_n}(u_B-u_S - u_B u_S, v_B-v_S)_{\Gamma(t)}\nonumber,\\
    S_h^n(u,v) &= \int_{I_n}s_h^n(t, u, v),\\
\begin{split}
    L_{h}^n(v) &= (u_{B,h}^-, v_B)_{\Omega(t_{n-1})} + (u_{S,h}^-, v_S)_{\Gamma(t_{n-1})} + \int_{I_n}(f_B,v_B)_{\Omega(t)}\\
    &+ \int_{I_n}(f_S,v_S)_{\Gamma(t)}.
\end{split}
\end{align}
For a derivation of this method, combine the derivation in Section~\ref{subsection:conservative-scheme} with the derivation of the method in~\cite{frachon2022cut} where Reynold's transport theorem is used for the surface PDE.\@ Since we want to develop a method that works also for higher order elements than linear, the ghost penalty stabilization alone is not enough to stabilize the discretization of the surface PDE~\cite{surface_stabilization}. For this coupled bulk-surface problem we combine the stabilization terms for the bulk problem with the stabilization terms proposed in~\cite{surface_stabilization} for the surface problem and get:
\begin{equation}
    \begin{aligned}
        s_h^n(t,u,v) &= \sum_{i=1}^m\Big(\sum_{F\in \F_h^n} \tau_F^i h^{2i-1} (\llbracket D_{\bm{n}_F}^i u_B\rrbracket_F, \llbracket D_{\bm{n}_F}^i v_B \rrbracket_F)_F\\
        &+  \sum_{F\in \F_{h,\Gamma}^n} \tau_{F,\Gamma}^i h^{2i-2} (\llbracket D_{\bm{n}_F}^i u_S\rrbracket_F, \llbracket D_{\bm{n}_F}^i v_S \rrbracket_F)\\
        &+ \tau_{\Gamma}^i h^{2i-2} (D^i_{\bm{n}_{\Gamma}} u_S, D^i_{\bm{n}_{\Gamma}}v_S)_{\Gamma(t)} \Big),
    \end{aligned}
\end{equation}
where $\tau_F^i, \tau_{F,\Gamma}^i, \tau_{\Gamma}^i > 0$, $\bm{n}_{\Gamma}$ denotes the unit normal vector on the interface, and $\F_{h,\Gamma}^n = \{F = K_1\cap K_2 : K_1, K_2 \in \T_{h,\Gamma}^n\}$. The corresponding patch-based stabilization is defined by 
\begin{equation}
    \begin{aligned}
    s_h^n(t,u,v) &= \sum_{F\in \F_h^n} \tau h^{-2}(\llbracket u_B\rrbracket_{\patch},\llbracket v_B\rrbracket_{\patch})_{\patch}+ \sum_{F\in \F_{h,\Gamma}^n} \tau_{\Gamma} h^{-3}(\llbracket u_S\rrbracket_{\patch},\llbracket v_S\rrbracket_{\patch})_{\patch}\\
    &+ \sum_{i=1}^m\tau_{\Gamma}^i h^{2i-2} (D^i_{\bm{n}_{\Gamma}} u_S, D^i_{\bm{n}_{\Gamma}}v_S)_{\Gamma(t)},
    \end{aligned}
\end{equation}
where $\tau, \tau_{\Gamma}, \tau_{\Gamma}^i > 0$ are given parameters.

In the case of macroelement stabilization, the set $\F_h^n$ is replaced with the set containing the internal edges of the macroelements as in Section~\ref{section:stabilization}. The construction of the macroelement partition of the active mesh $\T_{h,\Gamma}^n$ (see~\eqref{eq:active-mesh-time-surf}) associated with the surface problem is similar to the bulk case but an element is classified as large if
\begin{equation}\label{eq:largeelsurf}
    \delta_s \leq \frac{|K\cap \Gamma(t_q)|}{h^{d_i}},  \quad \text{for all } t_q\in Q_n. 
\end{equation}
Here, $d_i$ is the dimension of the surface $\Gamma$. In the model in two space dimensions $d_i=1$. As in Section~\ref{section:stabilization} a macroelement partition of $\T_{h,\Gamma}^n$ can be constructed and the set of stabilized edges $\F_{h,\Gamma}^n$ consist of internal faces of macroelements and stabilization is never applied between two macroelements.

We avoid selecting $\delta_s>0$ too large to prevent situations where no cut elements are classified as large. In the worst case, one can always apply full stabilization or connect as many small elements such that 
\begin{equation}
    \delta_s \leq \frac{|M\cap \Gamma(t_q)|}{h^{d_i}},  \quad \text{for all } t_q\in Q_n. 
\end{equation}

\subsection{Non-conservative scheme}

The non-conservative scheme from~\cite{hansbo2016cut} that we compare with reads: Given $u_h^-=u_h(t_{n-1}^-, \bm{x})$ find $u_h\in W_{h,k,m}^n$ such that
\begin{equation}
    A^n(u_h,v_h)+S^n(u_h,v_h) = L^n(v_h),\quad \text{for all }v_h\in W_{h,k,m}^n,
\end{equation}
where
\begin{align}
\begin{split}
    A^n(u, v) &= (u_B, v_B)_{\Omega(t_{n-1})} + (u_S,v_S)_{\Gamma(t_{n-1})} + \int_{I_n}(\partial_t u_B,v_B)_{\Omega(t)}\\
    &+ \int_{I_n}(\partial_t u_S, v_S)_{\Gamma(t)}\label{eq:A-coupled-nc} + \int_{I_n}(\bm{\beta}\cdot \nabla u_B,  v_B)_{\Omega(t)}\\
    &+ \int_{I_n}(\bm{\beta}\cdot\nabla u_S, v_S)_{\Gamma(t)} + \int_{I_n}((\nabla_{\Gamma}\cdot\bm{\beta})u_S,v_S)_{\Gamma(t)} \\
    &+ \int_{I_n}(D\nabla u_B, \nabla v_B)_{\Omega(t)}
    + \int_{I_n}(D_{\Gamma}\nabla_{\Gamma} u_S, \nabla_{\Gamma} v_S)_{\Gamma(t)}\\ 
    &+ \int_{I_n}(u_B-u_S-u_B u_S, v_B-v_S)_{\Gamma(t)},
\end{split}\\[1ex]
    S^n(u,v) &= \int_{I_n}s_h^n(t, u, v),\\
\begin{split}
    L^n(v) &= (u_{B,h}^-, v_B)_{\Omega(t_{n-1})} + (u_{S,h}^-, v_S)_{\Gamma(t_{n-1})} + \int_{I_n}(f_B,v_B)_{\Omega(t)}\\
    &+ \int_{I_n}(f_S,v_S)_{\Gamma(t)}.
\end{split}
\end{align}

Note that the last terms in~\eqref{eq:A-coupled-c} and~\eqref{eq:A-coupled-nc} are non-linear, hence Newton's method is used to solve the system of equations in each time step. Details can be found in e.g.~\cite{hansbo2016cut}.

\subsection{Numerical example}\label{subsection:coupled-numerical-example}

This example is a variant of Example 1 in~\cite{hansbo2016cut}. The boundary evolves with the same velocity field as in Example~\ref{subsection:ex1}, the diffusion coefficients are chosen as $D=0.01$, $D_{\Gamma}=1$. The force functions are computed such that the exact solutions are given by $u_B = 0.5 + 0.4\cos(\pi x)\cos(\pi y)\cos(2\pi t)$, and $u_S = (u_B + \bm{n}\cdot D \nabla u_B)/(1+u_B)$, where $\bm{n} = \begin{pmatrix} x-x_c(t), &y-y_c(t)
\end{pmatrix}^T/\sqrt{(x-x_c(t))^2+(y-y_c(t))^2}$ and $x_c$ and $y_c$ are given as in Example~\ref{subsection:ex1}. For all studies in this example, we set $\Delta t_n = h/4$.

\begin{figure}
    \centering
    \includegraphics[width=0.4\linewidth]{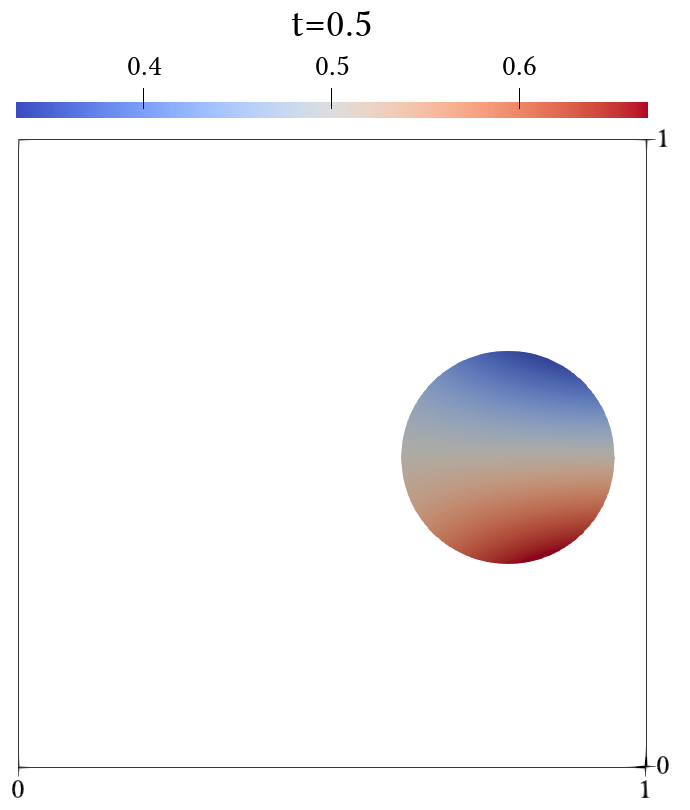}\hspace{15mm}
    \includegraphics[width=0.4\linewidth]{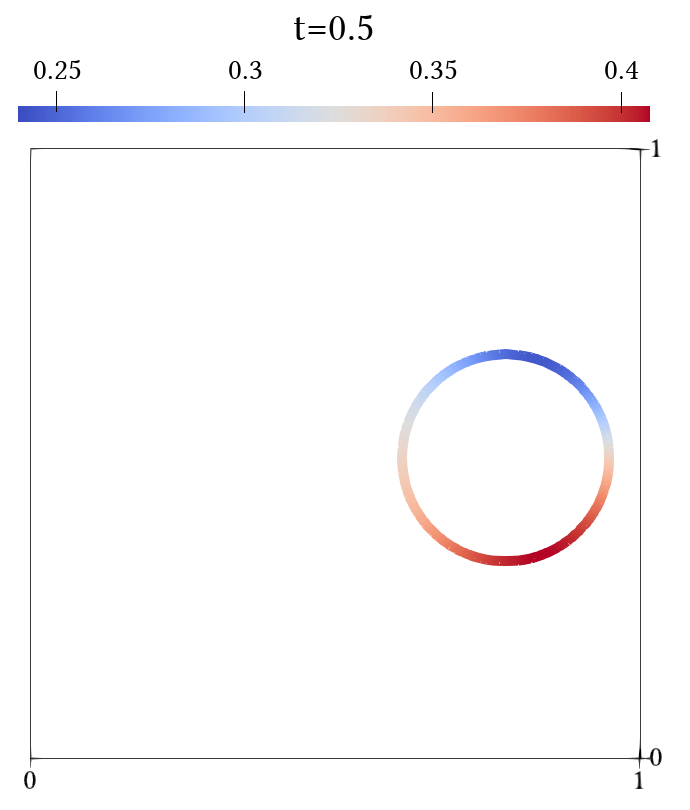}%
    \caption{Example~\ref{subsection:coupled-numerical-example}: Numerical solution in $t=0.5$. The bulk solution is shown to the left, and the surface solution to the right.}\label{fig:coupled-sol}
\end{figure}

\begin{figure}
    \centering
    \includegraphics[width=1\linewidth]{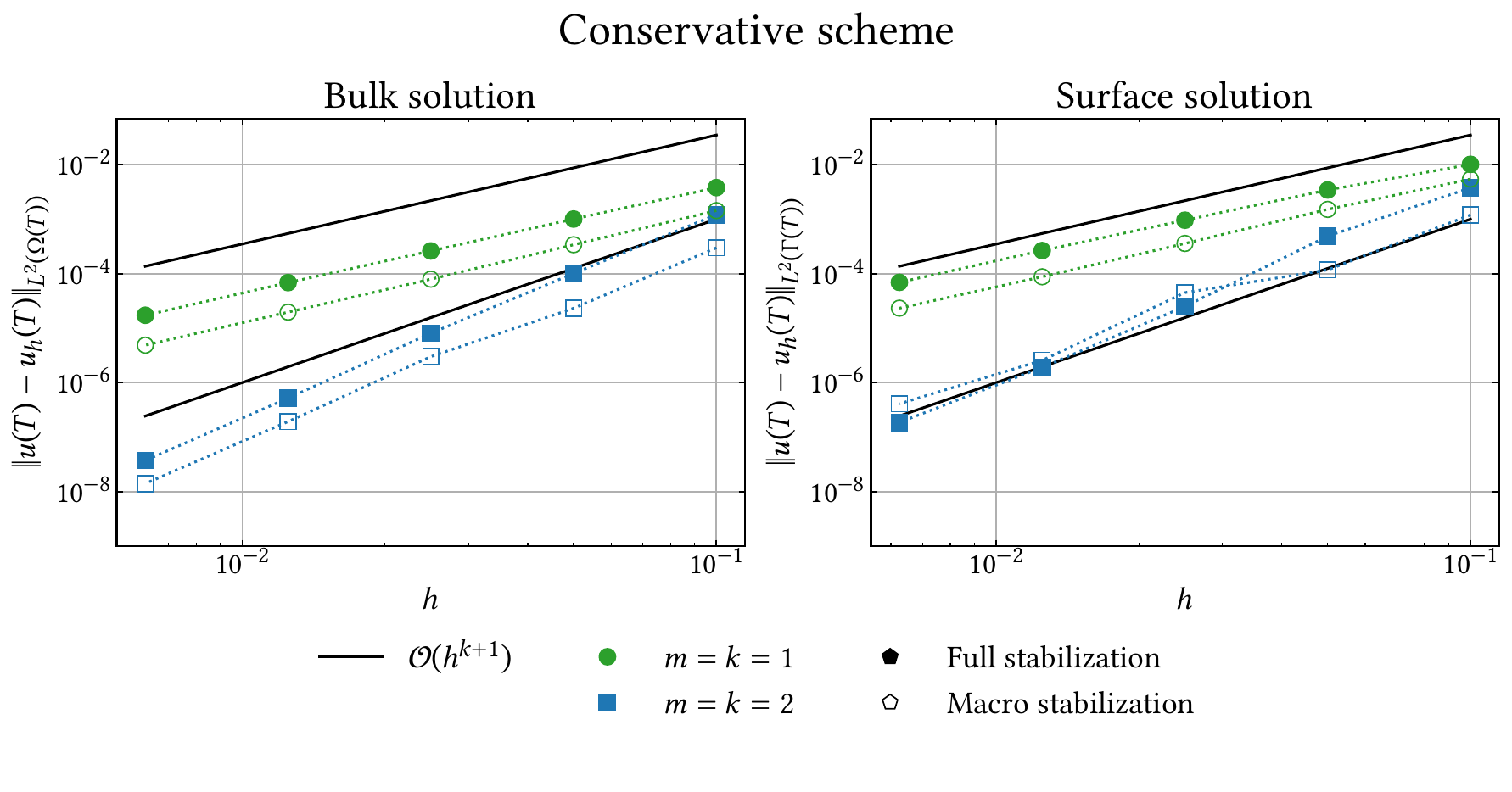}
    \caption{Example~\ref{subsection:coupled-numerical-example}: $L^2$ error for the bulk and surface solutions versus $h$ for the conservative scheme with $T=0.1$.}\label{fig:coupled_convergence_conservative}
\end{figure}

\begin{figure}
    \centering
    \includegraphics[width=1\linewidth]{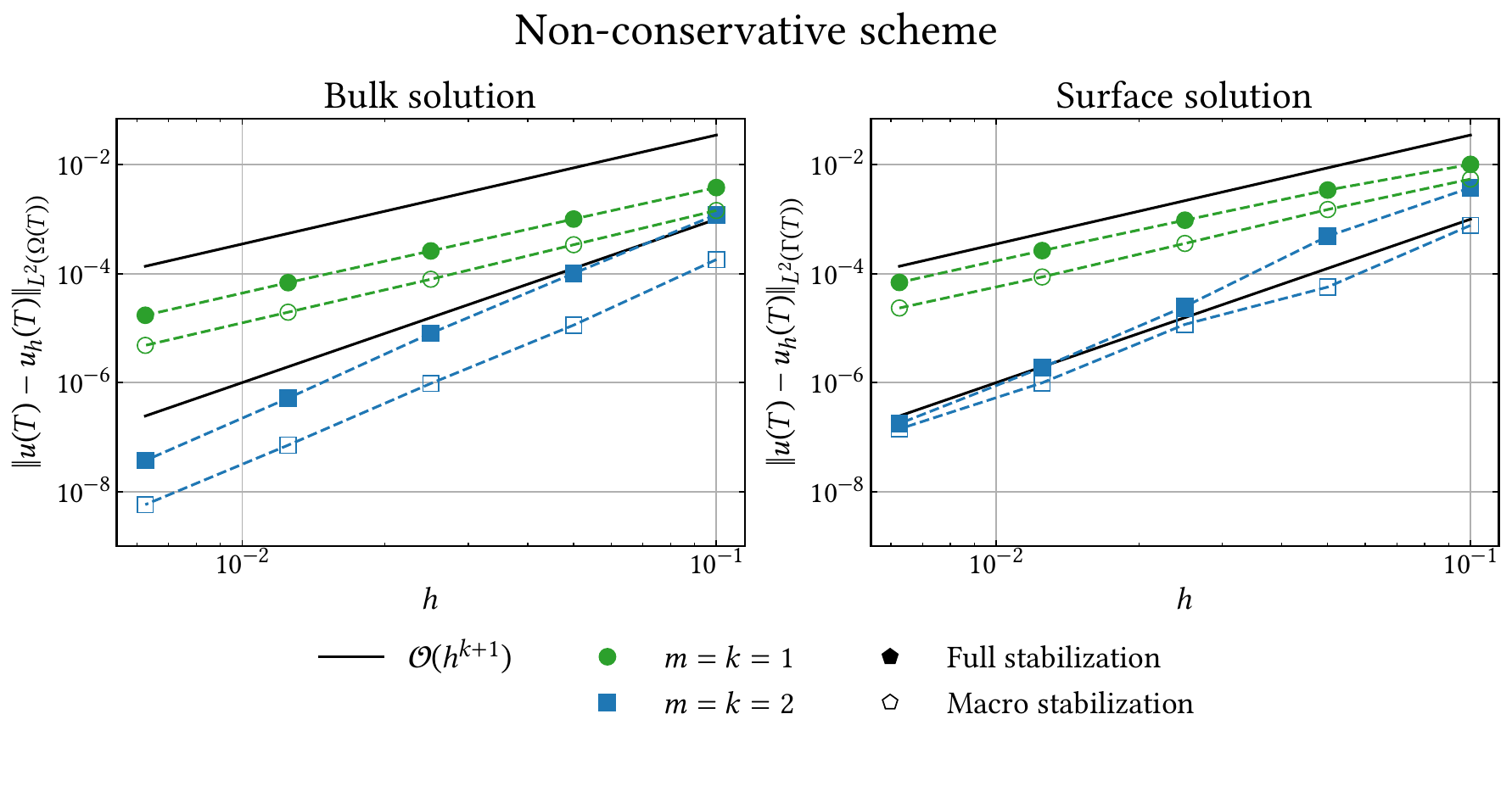}
    \caption{Example~\ref{subsection:coupled-numerical-example}: $L^2$ error for the bulk and surface solutions versus $h$ for the non-conservative scheme with $T=0.1$.}\label{fig:coupled_convergence_classical}
\end{figure}

\begin{figure}
    \centering
    \includegraphics[width=.7\linewidth]{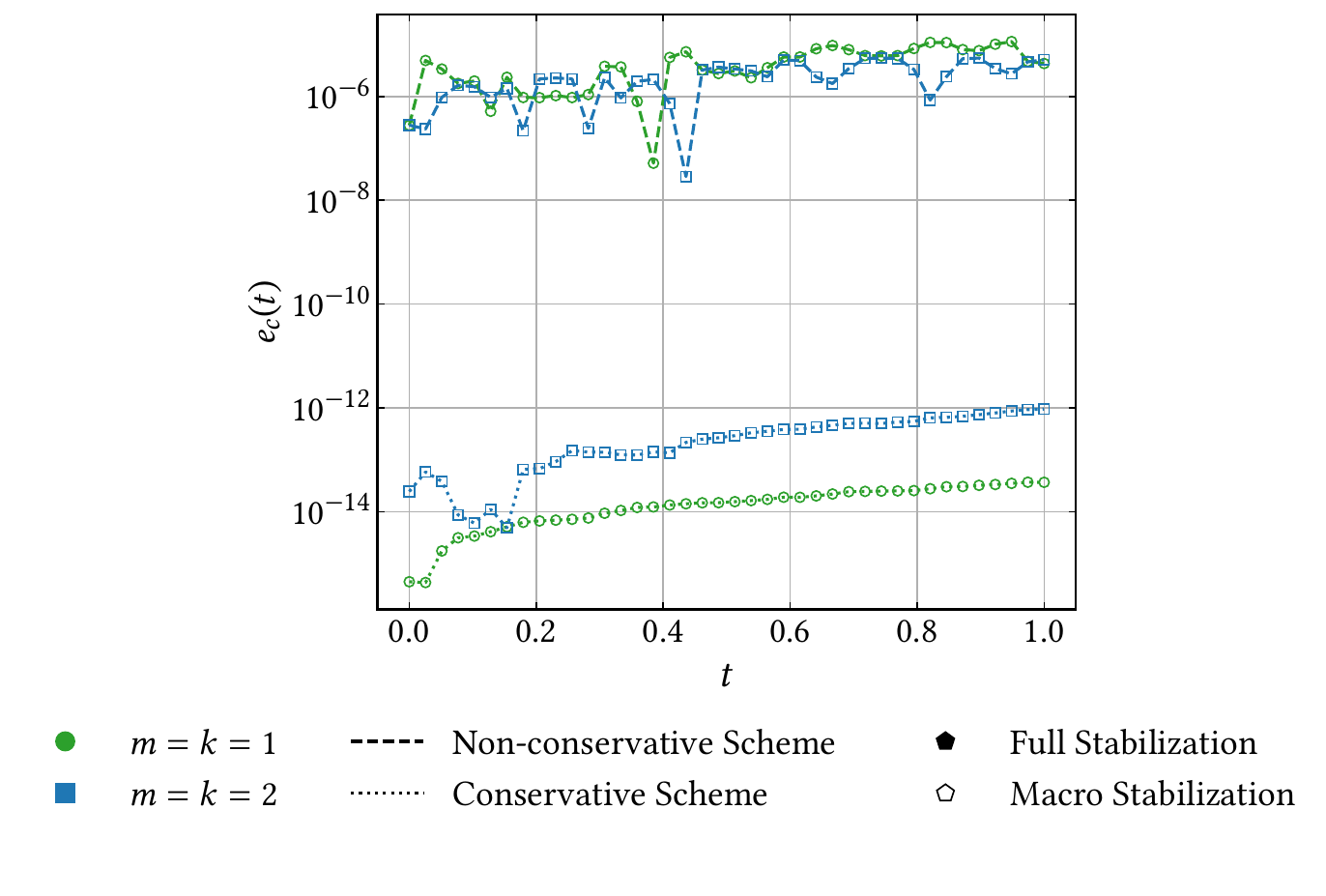}
    \caption{Example~\ref{subsection:coupled-numerical-example}: Conservation error $e_c(t)$ versus time $t$ with $h=0.1$, $T=1$.}\label{fig:coupled_conservation}
\end{figure}

The numerical solution is illustrated in Figure~\ref{fig:coupled-sol}. The convergence of the method is illustrated in Figure~\ref{fig:coupled_convergence_conservative} and Figure~\ref{fig:coupled_convergence_classical}. The conservation error is illustrated in Figure~\ref{fig:coupled_conservation}. The parameters are chosen as $\tau=\tau_{\Gamma}=\tau_{\Gamma}^i=1$ in all cases except for when using the conservative scheme with macroelement stabilization for $m=k=2$, then we use $\tau_{\Gamma}=10$ instead. The macroelement parameters are chosen as $\delta_B=0.7$ and $\delta_S=0.5$ for the bulk and surface respectively. When solving for the conservative scheme with $m=k=2$ and full stabilization, we use $N_t=20$. We observe optimal convergence for both schemes. The conservation error is of the order of machine precision for the conservative scheme also in this coupled problem.

\section{Conclusion}\label{section:conclusion}

Utilizing Reynold's transport theorem, we have presented an unfitted method that naturally conserves mass for the convection-diffusion equation in an evolving domain. The paper is focused on the problem posed in a bulk domain but also considers a non-linear coupled bulk-surface problem modeling the evolution of soluble surfactants. The method is shown to be optimal in terms of convergence and conditioning, and the conservation error is shown to be of the order of machine epsilon.

The macroelement stabilization we propose stabilizes the method more efficiently since stabilization is applied only where needed. This allows us to ensure stability for the conservative method with high-order elements independently of cut configurations. It yields control of the condition number, makes the scheme more robust with respect to large stabilization constants, and increases the sparsity of the matrix since stabilization is not applied between macroelements.

In this work, the velocity field advecting the surfactant and boundary was assumed to be known. In the future, we aim at coupling the presented method with a pointwise divergence-free unfitted method for the Stokes equations~\cite{frachon2023divergencefree}.

\section*{Acknowledgments}
This research was supported by the Swedish Research Council Grant No. 2018--04192, 2022--04808, and the Wallenberg Academy Fellowship KAW 2019.0190. The authors would like to thank Thomas Frachon for his help with some implementational aspects of the method.

\bibliographystyle{elsarticle-num} 
\bibliography{references.bib}

\end{document}